\documentclass[10pt,a4paper]{amsart} 
\usepackage{amsfonts,amssymb,amsthm,amsmath,ifpdf,stmaryrd,amsrefs,dsfont,esint,mathrsfs,dsfont}
\usepackage[hidelinks]{hyperref}
\usepackage{enumitem} 

\usepackage{tikz}
\usetikzlibrary{matrix,arrows}
\usepackage{tikz-cd}

\usepackage[overload]{textcase}

\ifpdf
  \usepackage[final,expansion=alltext,protrusion=alltext]{microtype} 
\fi

\usepackage{xargs} 
\usepackage{mathtools} 
\usepackage{ifthen} 

\theoremstyle{plain}

\newtheorem*{Th*}{Theorem}
\newtheorem*{ThA}{Theorem A}
\newtheorem*{ThB}{Theorem B}

\newtheorem*{Cor*}{Corollary}

\theoremstyle{definition}

\theoremstyle{remark}

\newif\ifTNS 
\TNStrue

\def\printtheoremname#1{\csname#1name\endcsname}
\def\printtheoremnames#1{\csname#1names\endcsname}

\def\thmref#1#2{\printtheoremname{#1}\ifTNS~\fi\ref{#1:#2}}

\def\uc#1#2{\MakeUppercase{#1}{#2}} 

\newcommand{\DefTheorem}[2]{\newenvironmentx{#1}[2][1=\empty,2=\empty]{%
    \ignorespaces%
    \ifx##2\empty%
      \begin{#2}%
    \else%
      \begin{#2}[{\uc##2}]%
    \fi%
    \ifx##1\empty%
      {}%
    \else%
      \label{#1:##1}%
    \fi%
    \ignorespaces}{\end{#2}\ignorespacesafterend}}

\newcommand{\prfof}[2]{\protect{Proof of~\thmref{#1}{#2}}}

\makeatother

\DefTheorem{Th}{theorem}
\DefTheorem{Prop}{proposition}
\DefTheorem{Cor}{corollary}
\DefTheorem{Lem}{lemma}
\DefTheorem{Def}{definition}
\DefTheorem{Rem}{remark}
\DefTheorem{Par}{para}
\DefTheorem{Not}{notation}
\DefTheorem{Exer}{exercise}
\DefTheorem{Ex}{example}
\DefTheorem{Cons}{construction}
\DefTheorem{Scho}{scholium}

\newenvironment{Par*}{\ignorespaces\noindent\ignorespaces}{\ignorespacesafterend}

\numberwithin{equation}{section}

\newcommand\Define[2][\empty]{\ignorespaces%
  \emph{#2}}%

\tikzset{
  commutative diagrams/.cd,
  arrow style=tikz,
  diagrams={>=stealth},
  shift up/.style={
    to path={([yshift=#1]\tikztostart.east) -- ([yshift=#1]\tikztotarget.west) \tikztonodes}},
  shift up left/.style={
    to path={([yshift=#1]\tikztostart.west) -- ([yshift=#1]\tikztotarget.east) \tikztonodes}},
  mathdouble/.style={-,double equal sign distance}
}

  {\ignorespaces\end{tikzpicture}%
  \end{minipage}\\%
  \ignorespacesafterend}

\def\ger{\mathfrak}

\newcommand\CategoryTypeface{\mathbf}
\def\cat{\CategoryTypeface}

\newcommand\SheafTypeface{\mathcal}
\def\sh{\SheafTypeface}

\def\DMO{\DeclareMathOperator}

\newcommand\ev{{\bar 0}}
\newcommand\odd{{\bar 1}}

\newcommand{\defi}{\coloneqq}     

\def\diff#1^#2{\ensuremath{\partial_{#1}^{#2}}}

\def\der#1/#2{\ifthenelse {\equal{#1}{}}
              {\ensuremath{\partial_{#2}{#1}}}
              {\ensuremath{\frac{\partial #1}{\partial #2}}}
        }

\def\derf#1/#2{\ifthenelse  {\equal{#1}{}}
              {\ensuremath{\frac{\partial #1}{\partial #2}}}
              {\ensuremath{\partial_{#2}{#1}}}
        }

\newcommand{\fa}{for all }

\newcommand{\fs}{for some }
\newcommand\mathfa[1][{}]{\quad\text{\fa{#1} }}

\newcommand{\scth}{such that }

\newcommand\eg{\emph{e.g.}~}
\newcommand\ie{\emph{i.e.}~}

\newcommand\via{\emph{via}~}
\newcommand\loccit{\emph{loc.~cit.}}
\newcommand\opcit{\emph{op.~cit.}}

\newcommand\etc{\emph{etc.}}

\def\multi(#1,#2){\ifthenelse {\equal{#1}{0}}
                {{\mathbb Z}_2^{#2}}
                {\ifthenelse{\equal{#2}{0}}
                      {{\mathbb N}_0^{#1}}
                      {\ensuremath{{{\mathbb N}_0^{#1}\!\times\!{\mathbb Z}_2^{#2}}}}
                }
        }

\newcommand\vphi{\varphi}
\newcommand\vrho{\varrho}

\newcommand\vkappa{\varkappa}
\newcommand\eps{\varepsilon}
\newcommand\nats{\mathbb{N}}
\newcommand\ints{\mathbb{Z}}

\newcommand\reals{\mathbb{R}}
\newcommand\cplxs{\mathbb{C}}
\newcommand\knums{\mathbb K}

\def\aff{{\mathbb A}}

\newcommand\sle{\leqslant}
\newcommand\sge{\geqslant}

\DMO\dom{\mathrm{dom}}
\DMO\rk{\mathrm{rk}}
\DMO\Ad{\mathrm{Ad}}
\DMO\ad{\mathrm{ad}}
\DMO\GL{\mathrm{GL}}
\DMO\id{\mathrm{id}}
\DMO\sdim{\mathrm{sdim}}
\DMO\sgn{\mathrm{sgn}}
\DMO\re{\mathrm{Re}}
\DMO\coker{\mathrm{coker}}
\DMO\im{\mathrm{im}}
\DMO\coim{\mathrm{coim}}
\DMO\codim{\mathrm{codim}}
\DMO\supp{\mathrm{supp}}
\DMO\str{\mathrm{str}}
\DMO\tr{\mathrm{tr}}

\DMO\HC{\cat{HC}}
\DMO\CW{\cat{CW}}
\DMO\SMan{\cat{SMan}}
\DMO\Sets{\cat{Sets}}
\DMO\SSp{\cat{SSp}}
\DMO\SVec{\cat{SVec}}
\DMO\Shv{\cat{Shv}}

\DMO\WF{\mathrm{WF}}
\DMO\Op{\mathrm{Op}}

\newcommand{\lBr}{[\kern-.65ex[}
\newcommand{\rBr}{]\kern-.65ex]}

\makeatletter
\newcommand\Size[7][1]{
                                 \ifx#20%
                                        \def\r@l{}\def\r@m{}\def\r@r{}%
                                 \else%
                                    \ifx#21%
                                           \def\r@l{\bigl}\def\r@r{\bigr}\def\r@m{\bigm}%
                                    \else%
                                           \ifx#22%
                                                 \def\r@l{\Bigl}\def\r@r{\Bigr}\def\r@m{\Bigm}%
                                            \else%
                                                 \ifx#23%
                                                        \def\r@l{\biggl}\def\r@r{\biggr}\def\r@m{\biggm}%
                                                  \else
                                                        \ifx#24%
                                                        \def\r@l{\Biggl}\def\r@r{\Biggr}\def\r@m{\Biggm}%
                                                        \fi%
                                                  \fi%
                                            \fi%
                                      \fi%
                                 \fi%
                                 \ifx#10%
                                       \def\r@m{}%
                                 \fi%
                                 \r@l#3{#4}\r@m#5{#6}\r@r#7%
}%

\makeatother

\makeatletter
\def\Set@Scallop[#1]#2#3{{#1}\Parens{#2}{#3}}
\newcommand\DeclareScalableOperator[2]{%
  \expandafter\def\csname#1\endcsname{\@ifnextchar[{{#2}\Set@Scallop}{{#2}\Set@Scallop[{}]}}
}
\makeatother

\def\DSO{\DeclareScalableOperator}

\DSO{Presh}{\cat{Presh}}
\DSO{Sh}{\cat{Sh}}

\DSO{ShHom}{\sh Hom} 
\DSO{ShGHom}{\underline{\sh Hom}} 
\DSO{ShEnd}{\sh End} 
\DSO{ShGEnd}{\underline{\sh End}} 
\DSO{ShDer}{\sh Der} 
\DSO{ShGDer}{\underline{\sh Der}} 
\DSO{ShCliff}{\sh C\text{\emph{liff}}} 

\DSO{Ct}{\mathcal{C}} 
\DSO{Lp}{\mathbf{L}} 
\DSO{Hom}{\mathrm{Hom}} 
\DSO{End}{\mathrm{End}} 
\DSO{Aut}{\mathrm{Aut}} 
\DSO{Uenv}{\mathfrak{U}} 
\DSO{Cuenv}{\widehat{\mathfrak{U}}} 
\DSO{ABer}{\Abs0{\mathrm{Ber}}}
\DSO{Proj}{\mathrm{Proj}}
\DSO{Ber}{\mathrm{Ber}}
\DSO{Vol}{\mathrm{Vol}}
\DSO{Form}{\Omega}
\DSO{Ind}{\mathrm{Ind}}
\DSO{Coind}{\mathrm{Coind}}
\DSO{Der}{\mathrm{Der}}
\DSO{GDer}{\underline{\mathrm{Der}}}
\DSO{Alg}{\mathrm{Alg}}
\DSO{GHom}{\underline{\mathrm{Hom}}} 
\DSO{GEnd}{\underline{\mathrm{End}}} 
\DSO{Maps}{\mathrm{Maps}}
\DSO{Cliff}{\mathrm{Cliff}} 
\DSO{Jac}{\mathrm{Jac}} 
\DSO{Dist}{\mathcal{D}'}
\DSO{Sw}{\mathscr{S}}
\DSO{TDi}{\mathscr{S}'}

\newcommand\Set[3]{
                                 \Size{#1}{\{}{#2}{|}{#3}{\}}%
}%
\newcommand\Dual[3]{
                                 \Size[0]{#1}{\langle}{#2}{,}{#3}{\rangle}%
}%
\newcommand\Parens[2]{
  \Size[0]{#1}{(}{#2}{}{}{)}
}

\newcommand\Braces[2]{
  \Size[0]{#1}{\{}{#2}{}{}{\}}
}

\newcommand\Norm[2]{
  \Size[0]{#1}{\lVert}{#2}{}{}{\rVert}
}
\newcommand\Abs[2]{
  \Size[0]{#1}{\lvert}{#2}{}{}{\rvert}
}
\newcommand\Span[2]{
  \Size[0]{#1}{\langle}{#2}{}{}{\rangle}
}
\makeatletter

\makeatletter

\newif\if@smallmat
\newif\if@none
\newif\if@paren
\newif\if@brack
\newif\if@brace
\newif\if@vline
                                 {\ifx#20%
                                        \@smallmattrue%
                                  \else%
                                         \@smallmatfalse
                                  \fi%
                                  \ifx#11%
                                         \@nonefalse\@parentrue\@brackfalse\@bracefalse\@vlinefalse%
                                  \else%
                                       \ifx#12%
                                            \@nonefalse\@parenfalse\@bracktrue\@bracefalse\@vlinefalse%
                                        \else%
                                            \ifx#13%
                                                 \@nonefalse\@parenfalse\@brackfalse\@bracetrue\@vlinefalse%
                                            \else%
                                                 \ifx#14%
                                                       \@nonefalse\@parenfalse\@brackfalse\@bracefalse\@vlinetrue
                                                 \else%
                                                       \ifx#15%
                                                             \@nonefalse\@parenfalse\@brackfalse\@bracefalse\@vlinefalse%
                                                       \else%
                                                             \@nonetrue\@parenfalse\@brackfalse\@bracefalse\@vlinefalse%
                                                       \fi%
                                                 \fi%
                                            \fi%
                                        \fi%
                                   \fi%
                                   \if@smallmat%
                                        \if@none%
                                             \begin{smallmatrix}%
                                        \else%
                                            \if@paren%
                                                  \left(\begin{smallmatrix}%
                                            \else%
                                                  \if@brack%
                                                          \left[\begin{smallmatrix}%
                                                  \else%
                                                          \if@brace%
                                                               \left\{\begin{smallmatrix}%
                                                          \else%
                                                               \if@vline%
                                                                    \left\lvert\begin{smallmatrix}%
                                                                \else%
                                                                    \left\lVert\begin{smallmatrix}%
                                                                \fi%
                                                          \fi%
                                                  \fi%
                                            \fi%
                                        \fi%
                                   \else%
                                        \if@none%
                                             \begin{matrix}%
                                        \else%
                                            \if@paren%
                                                  \begin{pmatrix}%
                                            \else%
                                                  \if@brack%
                                                          \begin{bmatrix}%
                                                  \else%
                                                          \if@brace%
                                                               \begin{Bmatrix}%
                                                          \else%
                                                               \if@vline%
                                                                    \begin{vmatrix}%
                                                                \else%
                                                                    \begin{Vmatrix}%
                                                                \fi%
                                                          \fi%
                                                  \fi%
                                            \fi%
                                        \fi%
                                   \fi}%
                                  {\if@smallmat%
                                        \if@none%
                                             \end{smallmatrix}%
                                        \else%
                                            \if@paren%
                                                  \end{smallmatrix}\right)%
                                            \else%
                                                  \if@brack%
                                                          \end{smallmatrix}\right]%
                                                  \else%
                                                          \if@brace%
                                                               \end{smallmatrix}\right\}%
                                                          \else%
                                                               \if@vline%
                                                                    \end{smallmatrix}\right\rvert%
                                                                \else%
                                                                    \end{smallmatrix}\right\rVert%
                                                                \fi%
                                                          \fi%
                                                  \fi%
                                            \fi%
                                         \fi%
                                   \else%
                                        \if@none%
                                             \end{matrix}%
                                        \else%
                                            \if@paren%
                                                  \end{pmatrix}%
                                            \else%
                                                  \if@brack%
                                                          \end{bmatrix}%
                                                  \else%
                                                          \if@brace%
                                                               \end{Bmatrix}%
                                                          \else%
                                                               \if@vline%
                                                                    \end{vmatrix}%
                                                                \else%
                                                                    \end{Vmatrix}%
                                                                \fi%
                                                          \fi%
                                                  \fi%
                                            \fi%
                                        \fi%
                                   \fi}%
                         
\makeatother

\def\clap#1{\hbox to 0pt{\hss#1\hss}} 

\begin{document}

\title[Harish-Chandra supermodules]{Fr\'echet globalisations of Harish-Chandra supermodules}

\author[Alldridge]
{Alexander Alldridge}

\address{Universit\"at zu K\"oln\\
Mathematisches Institut\\
Weyertal 86-90\\
50931 K\"oln\\
Germany}
\email{alldridg@math.uni-koeln.de}

\thanks{The author was supported by Deutsche Forschungsgemeinschaft, grant nos.~AL 698/3-1 (Heisenberg grant), ZI 513/2-1 (Leibniz prize to M.~Zirnbauer), SFB/TRR 12 ``Symmetries and University in Mesoscopic Systems'', SFB/TRR 183 ``Entangled States of Matter'', and by the Institutional Strategy of the University of Cologne within the German Excellence Initiative}



\begin{abstract}
  For any Lie supergroup whose underlying Lie group is reductive, we prove an extension of the Casselman--Wallach globalisation theorem: There is an equivalence between the category of Harish-Chandra modules and the category of $SF$-representations (smooth Fr\'echet representations of moderate growth) whose module of finite vectors is Harish-Chandra. As an application, we extend to Lie supergroups a general general form of the Gel\cprime fand--Kazhdan criterion due to Sun--Zhu. 
\end{abstract}

\maketitle

\section*{Introduction}

In the study of continuous representations of non-compact real-reductive Lie groups $G_0$, a fundamental obstacle is that almost all representations of interest are infinite-dimensional. A basic tool, which reduces many analytic questions to algebraic ones, is the passage to the module of $K_0$-finite vectors. The fundamental Casselman--Wallach theorem \cites{casselman,wallach2} guarantees that every Harish-Chandra $(\ger g_0,K_0)$-module occurs in this way. This is essential, in particular in applications to the  classification problem for irreducible unitary representations. 

Lie supergroups were introduced by Berezin, Kostant, and Le\u\i{}tes \cites{berezin-kac,berezin-leites,kostant} in the 1970s as a mathematical framework for the study of the supersymmetries occurring in quantum field theory. Lie superalgebras had entered the stage three decades earlier, in the work of Whitehead, to gain more prominence the work of Fr\"olicher--Nijenhuis on the deformation of complex structure, and subsequently in Gerstenhaber's work on the deformation of rings and algebras. Through the applications in physics, the subject of Lie superalgebras representations has come to the fore, and is at present well-established in both mathematics and physics, with a literature far too extensive to cite; compare the paper \cite{cns} for the development up to the 1970s, and the monographs \cites{cheng-wang,musson} for an up-to-date account of the theory. On the level of Lie supergroups, there is a sizeable literature in physics, but the subject has been hardly studied from a mathematical perspective. 

Most mathematical works (\eg Refs.~\cites{dz,furutsu-nishiyama,jakobsen,nishiyama,pj}) consider the unitarisable Harish-Chandra modules, without exploring the issue whether they arise as the space of finite vectors of some `global' representation. A `global' perspective was taken by Dobrev--Petkova \cite{dp-func}, who realise induced representations of the supergroup $\mathrm{SU}(2,2|N)$ on spaces of superfunctions. They classify unitary irreducible represenations of positive energy \cites{dp-pe1,dp-pe2}, thereby extending previous work of Flato--Fronsdal \cite{ff} for $N=1$. On the basis of these seminal ideas, Carmeli--Cassinelli--Toigo--Varadarajan \cite{cctv} introduce a notion of unitary Lie supergroup representations for arbitary Lie supergroups. This has spawned a flurry of further investigation \cites{salmasian-vn,ns-pd,ns-cone,mns}. 

Meanwhile, beyond the obvious fact that non-unitary representations may occur as intermediates in the study of unitary ones, it has become clear that unitary representations alone are insufficient for the purpose of Fourier--Plancherel decomposition, even in simple cases \cite{ahl}. This is confirmed by applications of supersymmetry to number theory and random matrices \cites{cfz,hpz}, as well as in physics, for instance in the study of the Chalker--Coddington model with point contacts \cite{bwz}.

Finally, as has become increasingly clear in recent investigations of the Gel\cprime fand and Gel\cprime fand--Kazhdan properties for pairs of Lie groups  beyond the setting of Riemannian symmetric pairs \cites{ag,ags,agrs,sun-zhu-annals} that Casselman--Wallach theory is eminently useful for the study of branching multiplicities. Here, we argue that similar statements hold true also for the setting of Lie supergroups.

Therefore, it seems paramount to study the globalisations of `algebraic' representations, irrespective of unitarity, for Lie supergroups. In this paper, we generalise the Casselman--Wallach theorem to Lie supergroups, as follows.

\begin{ThA}
  Let $G$ be a Lie supergroup whose underlying Lie group $G_0$ is almost connected and real reductive, $\ger g$ its Lie superalgebra, and let $K_0\subseteq G_0$ be maximally compact. Then any Harish-Chandra $(\ger g,K_0)$-module has a unique $SF$-globalisation. 

  This defines an additive equivalence between the category $\HC(\ger g,K_0)$ of Harish-Chandra $(\ger g,K_0)$-modules and the category $\CW(G)$ of $SF$-representations of $G$ whose module of $K_0$-finite vectors is Harish-Chandra.
\end{ThA}

Here, we follow Ref.~\cite{bk} in using the term `$SF$-re\-pre\-sen\-tation' (resp.~`$F$-re\-pre\-sen\-ta\-tion') instead of `smooth Fr\'echet representation of moderate growth' (resp.~`Fr\'echet representation of moderate growth').

As an application of our results on globalisation, we study the Gel\cprime fand--Kazhdan property for pairs of supergroups, to arrive by the following version of the Gel\cprime fand--Kazhdan criterion, which generalises that given recently by Sun--Zhu \cite{sun-zhu}.

\begin{ThB}
  Let $H_1,H_2\subseteq G$ be closed subsupergroups, $\chi_i$, $i=1,2$, characters of $H_i$, $i=1,2$, and $\sigma$ an antiautomorphism of $G$. Assume that any even relatively $(\chi_1^{-1}\otimes\chi_2^{-1})$-invariant tempered superfunction $G$ that is a joint eigenvector of all even $G$-invariant $D\in\Uenv0{\ger g}$ is fixed by $\sigma$.

  Then, for any contragredient pair $(E,F)$ of $F$-representations of $G$ \scth $E_\infty$ and $F_\infty$ are irreducible $G$-representations whose modules of $K_0$-finite vectors are Harish-Chandra, we have 
  \[
    \dim\Hom[_{H_1}]0{E_\infty,\chi_1}\dim\Hom[_{H_2}]0{F_\infty,\chi_2}\sle1.
  \]
\end{ThB}

\medskip\noindent
Theorem A (\thmref{Th}{cw}) is derived in the framework of convolution algebras of Schwartz functions, as used by Bernstein--Kr\"otz \cite{bk} for Lie groups in their proof of a Casselman--Wallach theorem for holomorphic families of Harish-Chandra modules. 

As it turns out, the framework of convolution superalgebras of superdistributions and Berezinian densities is well-adapted to the study of the classes of continuous and weakly smooth representations, introduced here. In fact, a version of the Dixmier--Malliavin theorem holds (\thmref{Prop}{dm}). 

Moreover, the convolution algebra of Schwartz--Berezin densities is equally well suited for the study of $F$- and $SF$-representations (or moderate growth representations) of Lie supergroups. Indeed, we prove a Schwartzian Dixmier--Malliavin theorem for $F$-representations (\thmref{Prop}{sw-dm}), generalising the corresponding result of Bernstein--Kr\"otz \cite{bk}. 

What makes the proof of our main results tick is the fact that all of the convolution superalgebras in question can be presented as coinduced modules (\thmref{Prop}{superdist}, \thmref{Cor}{ber-dens}, \thmref{Prop}{sw-iso}, and \thmref{Prop}{sw-conv}), allowing for a passage from Lie supergroups to supergroup pairs. We can thus reduce many analytic questions to the underlying Lie group and use Hopf algebraic methods of computation to arrive by our conclusions. 

The expression of the convolution superalgebras \via coinduced modules, whilst preserving the convolution product, is, however, a non-trivial fact. It is based on an extension of Bruhat's regularity theorem for left-invariant distributions (\thmref{Prop}{bruhat}), which, together with dualising module techniques, implies an expression of the invariant Berezin density in terms of the Haar density on the underlying Lie group (\thmref{Prop}{invber}). Such an expression was previously only known in very special cases, where, in particular, the `odd modular function' is trivial \cite{cz}. The present result is far more general and covers all kinds of Lie supergroups, including non-basic classical and even non-simple cases.

In the final Section \ref{sec:gk}, we apply our results to the generalisation of the Gel\cprime fand--Kazhdan criterion in Theorem B (\thmref{Th}{super-gk}). The setting of Sun--Zhu \cite{sun-zhu} goes over more or less verbatim, due to our extension of the Casselman--Wallach theory. 

We do not yet view these last results as the definitive statements on multiplicity freeness for supergroups. Indeed, whereas we have focused here on the extension of phenomena from the purely even setting, there are many issues special to the super case that yet need to be addressed, such as $Q$ type modules and lack of semi-simplicity at the level of finite-dimensional modules. Moreover, non-trivial examples that verify the assumptions of Theorem B have yet to be supplied, and we intend to study this question in future work. However, the ease with which at least the purely even results transfer to the super case is to our mind a strong indication to the utility of the super Casselman--Wallach Theorem A. 

\medskip\noindent
\emph{Acknowledgements.} We extend our deep gratitude to the Institute for Theoretical Physics at the University of Cologne for its hospitality during the preparation of this article. We wish to thank the organisers of the Oberwolfach conference on ``Representations of Lie groups and supergroups'', Joachim Hilgert, Toshiyuki Kobayashi, Karl-Hermann Neeb, and Tudor Ratiu, where the results presented in this article were first announced. Moreover, we wish to thank the anonymous referees for their helpful comments. 

\section{Supergroup representations}

In this section, we collect some preliminary material on supergroups and their representations. 

\subsection{Preliminaries and notation}\label{subs:prelim}

Concerning supermanifolds, we will essentially work in the standard framework of Berezin and Le\u\i{}tes \cite{berezin-leites}, and use standard facts concerning it, to be found in Refs.~\cites{leites,manin,deligne-morgan,ccf}. We give some basic definitions to fix our terminology and to clarify in which points we deviate from this literature. 

We consider sheaves of Abelian groups and will denote them by calligraphic Roman letters $\sh E,\sh F,\sh O$, \etc{} The set of sections of a sheaf $\sh F$ will be denoted by $\Gamma(\sh F)$. The \Define{support} of a section $f\in\Gamma(\sh F)$ is $\Set0x{f_x\neq0}$ where $f_x$ is the germ at $x$. We denote by $\Gamma_K(\sh F)$ the set of sections with support contained in $K$, and by $\Gamma_c(\sh F)$ the set of sections with compact support. 

Let $\knums$ be the field $\reals$ of real or the field $\cplxs$ of complex numbers. Consider the category of $\knums$-superspaces: Its objects are pairs $X=(X_0,\sh O_X)$ comprised of a topological space $X_0$ and a sheaf $\sh O_X$ on $X_0$ of supercommutative $\knums$-superalgebras with local stalks; its morphisms $\vphi:X\longrightarrow Y$ are pairs $(\vphi_0,\vphi^\sharp)$ consisting of a continuous map $\vphi_0:X_0\longrightarrow Y_0$ and an even unital morphism of $\knums$-superalgebra sheaves $\vphi^\sharp:\sh O_Y\longrightarrow(\vphi_0)_*\sh O_X$ where $(\vphi_0)_*$ denotes the direct image functor. By a standard adjunction \cite{bredon}*{Chapter I.3}, we may equivalently consider $\vphi^\sharp$ as a morphism of sheaves $\vphi_0^{-1}\sh O_Y\longrightarrow\sh O_X$, where $\vphi_0^{-1}$ is the inverse image functor. Given some finite-dimensional super-vector space $V=V_\ev\oplus V_\odd$ over $\reals$, together with a compatible $\knums$-structure on the odd part $V_\odd$, we define the \Define{affine superspace} $\aff(V)$ by
\[
  \aff(V)_0\defi V_\ev,\quad\sh O_{\aff(V)}\defi\sh C^\infty_{V_\ev}\otimes_\reals\textstyle\bigwedge_\knums(V_\odd^*).
\]
Here, $\sh C^\infty_{V_\ev}$ denotes the sheaf of smooth real-valued functions on $V_\ev$, and $\bigwedge_\knums(V_\odd^*)$ denotes the exterior algebra of the $\knums$-vector space $V_\odd^*$. Here and in what follows, we denote the homogeneous parts of a given grading over $\ints/2\ints=\{\ev,\odd\}$ by the subscripts ${}_\ev$ (even) and ${}_\odd$ (odd).

Given a $\knums$-superspace $X$, an \Define{open subspace} is one of the form $X|_U\defi(U,\sh O_X|_U)$ for some open subset $U\subseteq X_0$. A $\knums$-superspace $X$ is called a \Define{supermanifold} if $X_0$ is Hausdorff and admits an open cover $(U_i)$ \scth for every index $i$, $X|_{U_i}$ is isomorphic to an open subspace of some affine superspace $\aff(V)$. In this case, for $x\in U_i$, the tuple $\dim_\reals V_\ev|\dim_\knums V_\odd$ is denoted by $\dim_xX$ and called the \Define{superdimension} of $X$ at $x$.

For $\knums=\reals$, one customarily calls supermanifolds as defined above \Define{real supermanifolds}; in the case $\knums=\cplxs$, they are called \Define{\emph{cs} manifolds} \cite{deligne-morgan}*{\S 4.8}. Notice that due to our preference for representations on complex vector spaces, we have a natural bias towards working in the latter setting. This is the main point in which we do not follow the standard texts. Most aspects of real supermanifolds carry over to the \emph{cs} case, with some notable exceptions related to real structures and the representability of vector bundles. We will take care to point these out to the reader.

\subsection{Supergroups and supergroup pairs}

It is known that the category of supermanifolds admits finite products \cite{leites}*{3.1.6}. Thus, group objects and their morphisms in this category are well-defined \cite{maclane}*{Chapter III.6}. A group object in the category of supermanifolds will be called a \Define{Lie supergroup} or simply a \Define{supergroup}. For $\knums=\reals$, these are real Lie supergroups, while for $\knums=\cplxs$, they are \emph{cs} Lie supergroups. 

\medskip\noindent
For applications to linear representations, the following definition proves useful.

\begin{Def}[sgrp-pair][supergroup pairs]
  Let $G_0$ be a real Lie group with Lie algebra $\ger g_0$, $\ger g$ be a Lie superalgebra over $\knums$ \scth $\ger g_\ev=\ger g_0\otimes_\reals\knums$, and $\Ad:G_0\longrightarrow\Aut0{\ger g}$ a smooth action of $G_0$ by Lie $\knums$-superalgebra automorphisms. We say that $(\ger g,G_0)$ (where the action is understood) is a \Define{supergroup pair} if the differential $d\Ad$ of $\Ad$ is the restriction of the bracket $[\cdot,\cdot]$ of $\ger g$ to $\ger g_0\times\ger g$.

  A \Define{morphism of supergroup pairs} $(\ger g, G_0)\longrightarrow (\ger h,H_0)$ consists by definition of a morphism $\vphi_0:G_0\longrightarrow H_0$ of real Lie groups and a $\vphi_0$-equivariant Lie $\knums$-superalgebra morphism $d\vphi:\ger g\longrightarrow\ger h$ \scth $d\vphi_0=\smash{d\vphi|_{\ger g_0}}$.
\end{Def}

In the literature, supergroup pairs are referred to as Harish-Chandra pairs. Since to our knowledge, Harish-Chandra never worked on supergroups, we prefer to use a less colourful nomenclature. 

The following proposition is due to Kostant \cite{kostant} and Koszul \cite{koszul} in the case $\knums=\reals$; see \cite{ccf}*{Chapter 7} for a detailed exposition. The extension to the case of $\knums=\cplxs$ presents no difficulty.

\begin{Prop}
  Consider the functor that assigns to a Lie supergroup $G$ the supergroup pair $(\ger g,G_0)$, where $G_0$ is the underlying Lie group of $G$, $\ger g$ is its Lie superalgebra, and $G_0$ acts on $\ger g$ by the natural adjoint action. 

  This functor defines an equivalence of the category of Lie supergroups and their morphisms with the category of supergroup pairs and their morphisms. 
\end{Prop}

\begin{Rem}
  In particular, we may associate with any real Lie supergroup $G$ the \emph{cs} Lie supergroup whose supergroup pair is $(G_0,\ger g\otimes_\reals\cplxs)$. On the level of superspaces, this sends $G$ to the complex superspace $(G_0,\sh O_G\otimes_\reals\cplxs)$.

  We are mainly interested in complex representations, so we consider the case of \emph{cs} Lie supergroups to be more relevant than the case of real Lie supergroups. Compare \cite{deligne-morgan}*{\S~4.9} for a list of five exemplary situations where it is more natural or even required to consider \emph{cs} manifolds instead of real supermanifolds. In particular, Example 4.9.3 (\opcit) describes a \emph{cs} Lie supergroup which does not admit a real form. By contrast, any complex Lie supergroup whose underlying Lie group has a real form, has a \emph{cs} form. (NB: It is known that there are nilpotent step-$2$ Lie algebras without a real form.)
\end{Rem}

\subsection{Smooth and continuous supergroup representations}

In what follows, let $G_0$ be a Lie group with Lie algebra $\ger g_0$. To fix our terminology, we recall the following somewhat standard definitions. 

\begin{Def}[smoothcontrep-group][continuous and smooth representations]
  Let $G_0$ be a Lie group, $E$ a topological vector space over $\knums$ and $G_0\times E\longrightarrow E$ a linear left action of $G_0$ on $E$. If the action is a continuous map, then we say that the induced map $\pi_0:G\longrightarrow\GL(E)$ is a \Define{continuous representation} of $G_0$ on $E$. 

  Let the topology on $E$ be locally convex. A vector $v\in E$ is called \Define{smooth} if the orbit map $\gamma_v:G_0\longrightarrow E:g\longmapsto\pi_0(g)v$ is a smooth map. For $x\in\ger g_0$, one defines
  \[
    d\pi_0(x)v\defi\frac d{dt}\Big|_{t=0}\pi_0(\exp(tx))v.
  \] 
  This defines an action of $\ger g_0$ on the space $E_\infty$ of all smooth vectors. One endows $E_\infty$ with the coarsest locally convex topology \scth \fa $u\in\Uenv0{\ger g_0}$, the linear map
  \[
    d\pi_0(u):E_\infty\longrightarrow E
  \]
  is continuous. Compare \citelist{\cite{bk}*{2.4.2} \cite{bruhat}*{\S 2.3} \cite{taylor-ncharman}*{Chapter 0} \cite{warner-harman1}*{Section 4.4.1}} for alternative definitions of the topology. The representation $\pi_0$ is called \Define{weakly smooth} if the canonical inclusion $E_\infty\longrightarrow E$ is an isomorphism of topological vector spaces. (For the reasons explained in \cite{bk}*{Remark 2.12}, we reserve the term \emph{smooth} for $F$-representations, to be defined below.)
\end{Def}

In what follows, let $G$ be a Lie supergroup with underlying Lie group $G_0$ and Lie superalgebra $\ger g$. We continue to denote the Lie algebra of $G_0$ by $\ger g_0$; in particular, $\ger g_0$ is a real form of $\ger g_\ev$, that is $\ger g_\ev=\ger g_0\otimes_\reals\knums$. We intend to give a definition of what a representation of $G$ is. For finite-dimensional representations, there are two possible ways to do this: A `piecemeal' definition in terms of pairs of representations of $G_0$ and $\ger g$ with a suitable compatibility, and a definition in functorial terms. 

To state this precisely, let $E$ be a super-vector space over $\knums$ of dimension $p|q$. Choose any homogeneous $\knums$-basis $(z_j)$ of $E$, and let $(z^j)$ be the dual basis. We let $E_\reals$ be $E$, where we forget the $\knums$-structure on $E_\ev$ and retain only the $\reals$-structure. If $\Abs0{z^j}=\ev$, then we decompose $z^j=x^j+iy^j$ where $x^j,y^j\in\Hom[_\reals]0{E_\ev,\reals}$. 

By \cite{leites}*{Theorem 2.1.7}, there is for any supermanifold $S$ a natural bijection
\begin{gather*}
  \Braces1{v:S\longrightarrow\aff(E_\reals)}\longrightarrow\Gamma(\sh O_{S,\ev,\reals}^{np})\times\Gamma(\sh O_{S,\odd}^q)\\
  v\longmapsto\Parens1{v^\sharp(x^1),v^\sharp(y^1),\dotsc,v^\sharp(x^p),v^\sharp(y^p),v^\sharp(z^{p+1}),\dotsc,v^\sharp(z^{p+q})}.  
\end{gather*}
Here, $n=\dim_\reals\knums$, it is understood that $y^j=0$ for $\knums=\reals$, and $\sh O_{S,\ev,\reals}$ denotes the subsheaf of \emph{real-valued} sections of $\sh O_{S,\ev}$. That is, the canonical image of any germ $f_s$ of a section in the residue field $\vkappa(s)=\sh O_{S,s}/\ger m_{S,s}=\cplxs$ is required to lie in the subfield $\reals$ at every point $s$. The above bijection is natural in the choice of bases.

For any supermanifold $S$, we define $\GL(E)(S)$ to be set of invertible even matrices $g=(g_{k\ell})$ with entries in $\Gamma(\sh O_S)$. There is an obvious way to define $\GL(E)$ on morphisms $S'\longrightarrow S$, turning the assignment $S\longmapsto\GL(E)(S)$ into a set-valued cofunctor on $\knums$-supermanifolds. Naturally, $\GL(E)$ is a functor in groups. 

The following proposition is a minor variation on a well-known statement and can be easily derived from standard facts on supergroup actions given ample exposition in the literature, see \eg \cite{ccf}*{Chapter 8}. Compare also the recent work of Ostermayr \cite{ostermayr}*{Section 2}.

\begin{Prop}[rep-equiv]
   Then the following data are in one-to-one correspondence:
  \begin{enumerate}
    \item\label{item:rep-quiv-i} Pairs $(d\pi,\pi_0)$ of graded linear representations $\pi_0$ of $G_0$ on $E$ and $G_0$-equivariant Lie superalgebra actions $\pi$ of $\ger g$ on $E$ with $d\pi_0=d\pi|_{\ger g_0}$;
    \item\label{item:rep-quiv-ii} left actions $a:G\times E\longrightarrow E$ that are $\knums$-linear over $G$, in the sense that there is some $g=(g_{k\ell})\in\GL(E)(G)$ \scth 
    \[
      \sum_kg_{kj}z^k=
      \begin{cases}
        a^\sharp(x^j)+ia^\sharp(y^j), &\text{if }j\sle p,\\
        a^\sharp(z^j), &\text{if }j>p.
      \end{cases}
    \]
  \end{enumerate}
\end{Prop}

\begin{Rem}
  For $\knums=\reals$, $\GL(E)$ is represented by the real Lie supergroup $\GL(E,\reals)$ whose underlying supergroup pair is $(\GL(E_\ev,\reals)\times\GL(E_\odd,\reals),\ger{gl}(E,\reals))$. Thus, the data in \thmref{Prop}{rep-equiv} \eqref{item:rep-quiv-ii} are just morphisms of Lie supergroups $G\longrightarrow\GL(E,\reals)$.

  On the other hand, for $\knums=\cplxs$ (the case of \emph{cs} Lie supergroups), the proposition does not admit a statement in terms of supergroup homomorphisms. Indeed, in this case, $\GL(E)$ is not representable in the category of supermanifolds (\ie \emph{cs} manifolds). Instead, if we extend $\GL(E)$ to a suitable subcategory of the category of $\cplxs$-superspaces which contains complex supermanifolds as a (full) subcategory, then $\GL(E)$ can be seen to coincide on this subcategory with the point functor of the \emph{complex} Lie supergroup $\GL(E,\cplxs)$, which is not a \emph{cs} manifold for $E\neq0$. See Ref.~\cite{ostermayr}*{Appendix 7.1}.

  Another way to see that the data in \thmref{Prop}{rep-equiv} \eqref{item:rep-quiv-i} do not correspond to the $G$-points of a representable functor is to remark that the even part of the complex Lie super\-al\-ge\-bra $\ger{gl}(E,\cplxs)$ is not the complexification of the Lie algebra of $\GL(E_\ev,\cplxs)\times\GL(E_\odd,\cplxs)$, considered as a real Lie group, and thus, these do not form a supergroup pair.
\end{Rem}

On grounds of the above equivalence, we adopt the following terminology.

\begin{Def}[smoothcontrep-sgroup][continuous and smooth supergroup representations]
  Let $E$ be a locally convex super-vector space (\ie $E$ is a locally convex vector space with a grading that exhibits $E$ as a locally convex direct sum). Assume given a continuous representation $\pi_0$ of $G_0$ on $E_\ev$ and a Lie superalgebra representation $d\pi$ of $\ger g$ on $E_\infty$ \scth the map $\ger g\times E_\infty\longrightarrow E_\infty:(x,v)\longmapsto d\pi(x)v$ is continuous.

  We say that $(d\pi,\pi_0)$ is a \Define{continuous $G$-representation} if $d\pi$ is $G_0$-equivariant, \ie
  \[
    d\pi\Parens1{\Ad(g)(x)}=\pi_0(g)d\pi(x)\pi_0(g^{-1})
  \]
  \fa $x\in\ger g$ and $g\in G_0$, and $d\pi|_{\ger g_0}=d\pi_0$. If in addition, $E$ is weakly smooth as a $G_0$-representation, then we call $E$ a \Define{weakly smooth $G$-representation}.
\end{Def}

The definition given above for continuous supergroup representations is compatible with the corresponding ones given in the literature for the case of unitary representations \citelist{\cite{cctv}*{2.3} \cite{mns}*{Definition 4.1}}.

\section{Convolution superalgebras and representations}

In what follows, let $G$ be a Lie supergroup, where $G_0$ is assumed to be $\sigma$-compact. Let $\ger g$ be its Lie superalgebra. In this section, we introduce a convolution superalgebra of compactly supported Berezinian densities on $G$ and show that there is a one-to-one correspondence between its (non-degenerate) representations and the smooth representations of $G$. 

To that end, we will identify the sheaf of Berezinian densities of $G$ within the sheaf of superdistributions as the $\ger g$-module induced from the sheaf of densities on the underlying Lie group $G_0$. We begin by discussing superdistributions.

\subsection{Superdistributions}\label{subs:superdist}

In this section, we introduce superdistributions on $G$, and show how to express them in terms of the underlying Lie group. 

For any open $U\subseteq G_0$, we endow $\sh O_G(U)$ with the locally convex topology generated by the seminorms 
\[
  p_{u,v,K}(f)\defi\sup\nolimits_{x\in K}\Abs1{(L_uR_vf)(x)}
\]
where $K\subseteq U$ is compact and $u,v\in\Uenv0{\ger g}$. Here, $L$ and $R$, respectively, denote the left and right regular representation. It is known \citelist{\cite{koszul}*{Section 1} \cite{ccf}*{Proposition 7.4.13}} that there is an isomorphism 
\begin{center}
  \begin{tikzcd}
    \sh O_G(U)\rar{\phi}&\GHom[_{\ger g_\ev}]0{\Uenv0{\ger g},\sh C^\infty(U,\knums)}
  \end{tikzcd}
\end{center}
given by 
\[
  \phi(f)(u;x)\defi(-1)^{\Abs0f\Abs0u}(R_uf)(x)
  \]
\fa $f\in\sh O_G(U)$, $u\in\Uenv0{\ger g}$, and $x\in U$. Here, the action of $\ger g_\ev$ on $\sh C^\infty(U,\knums)$ is by left-invariant differential operators (\ie infinitesimal right translations), and the algebra product is expressed on the right-hand side by the rule
\[
  fh=m\circ(f\otimes h)\circ\Delta,
\]
where $m$ denotes multiplication in $\sh C^\infty$ and $\Delta$ denotes comultiplication in $\Uenv0{\ger g}$. For future use, we note that the multiplication morphism $m$ is given in terms of the isomorphism $\phi$ by 
\begin{equation}\label{eq:mult}
  \begin{aligned}[c]
    \phi(m^\sharp(f))(u\otimes v;g,h)&=\phi(f)\Parens1{\Ad(h^{-1})(u)v;gh}\\
   &=(-1)^{(\Abs0u+\Abs0v)\Abs0f}\Parens1{L_{S(\Ad(g)(u))}R_vf}(gh).
  \end{aligned}
\end{equation}

Since $\Uenv0{\ger g}=\Uenv0{\ger g_\ev}\otimes\bigwedge\ger g_\odd$ as graded $\ger g_\ev$-modules \cite{scheunert}*{I.2.3}, we have
\[
  \sh O_G(U)\cong\sh C^\infty(U,\knums)\otimes_\knums\textstyle\bigwedge(\ger g_\odd)^*.
\]
Since the Grassmann factor is finite-dimensional, one readily checks that is an isomorphism of locally convex super-vector spaces, where $\sh C^\infty(U,\knums)$ is given the usual topology of uniform convergence with all derivatives on compact subsets. In particular, $\sh O_G(U)$ is an $m$-convex Fr\'echet algebra \cite{michor-cinfty}*{2.2}. Here, we recall that a locally convex algebra is called $m$-convex if its topology is generated by a set of submultiplicative seminorms. 

Similarly, we give $\Gamma_c(\sh O_G)$ the locally convex inductive limit topology for the embeddings of the subspaces $\Gamma_K(\sh O_G)$ of sections $f$ with support $\supp f\subseteq K$, where $K\subseteq G_0$ is compact. (See Subsection \ref{subs:prelim} for the notation.) The latter are given the relative topology induced by $\Gamma(\sh O_G)$. Then $\Gamma_c(\sh O_G)$ is an LF space, and the multiplication is jointly continuous. 

\begin{Def}[dist][superdistributions]
  For any open $U\subseteq G_0$, define
  \[
    \sh Db_G(U)\defi\Gamma_c(\sh O_G|_U)',
  \]
  the strong continuous dual space. Since $\sh O_G$ is a $c$-soft sheaf, we have by \cite{bredon}*{Chapter V, \S 1, Proposition 1.6} that $U\longmapsto\Gamma_c(\sh O_G|_U)$ is a flabby cosheaf. The corestriction maps are continuous by the definition of the topology. Thus, $\sh Db_G$ is a sheaf of locally convex super-vector spaces, called the \Define{sheaf of superdistributions}. In particular, we let $\sh D'(G)\defi\Gamma(\sh Db_G)=\Gamma_c(\sh O_G)'$. 

  The sheaf $\sh Db_G$ is naturally a right $\sh O_G$-module by 
  \[
    \Dual0{\mu f}\vphi\defi\Dual0\mu{f\vphi},
  \]
  where $\Dual0\cdot\cdot$ denotes the canonical pairing between $\sh Db_G(U)$ and $\Gamma_c(\sh O_G|_U)$. 

  The Lie supergroup $G$ acts from the left on $\sh Db_G$, where the $G_0$- and $\ger g$-action are given respectively by
  \[
    \Dual0{L_g\mu}\vphi\defi\Dual0\mu{L_{g^{-1}}\vphi},\quad\Dual0{L_x\mu}\vphi\defi-\Dual0\mu{L_x\vphi}.
  \]
  Here, in terms of the isomorphism $\phi$, we have 
  \[
    \phi(L_g\vphi)(u;h)=\phi(\vphi)(u;g^{-1}h),\quad
    \phi(L_x\vphi)(u;h)=-\phi(\vphi)(\Ad(h^{-1})(x)u;h).
  \]
\end{Def}

In what follows, if $H$ is a subsupergroup of $G$ and $\sh A$ is a subalgebra of $\sh O_G$, we will call a sheaf on $G_0$ with a left $H$-action commuting with a right $\sh A$-action an $(H,\sh A)$-module. Thus, $\sh Db_G$ is a $(G,\sh O_G)$-module.

In the following proposition, recall that any supermanifold $X$ comes with a natural embedding of the underlying manifold $X_0$, denoted by $j_{X_0}:X_0\longrightarrow X$. The underlying map of $j_{X_0}$ is the identity; the sheaf map $\smash{j_{X_0}^\sharp}$ assigns to any superfunction $f$ its underlying function $f_0$.

\begin{Prop}[superdist]
  Let $\sh Db_{G_0}$ be the sheaf of superdistributions on $G_0$. There is an isomorphism of $(G,\sh O_{G_0})$-modules
  \begin{center}
    \begin{tikzcd}
      \Uenv0{\ger g}\otimes_{\Uenv0{\ger g_\ev}}\sh Db_{G_0}\rar{}&\sh Db_G,
    \end{tikzcd}
  \end{center}
  given by 
  \begin{equation}\label{eq:tensor-dist}
    \Dual0{u\otimes\mu}\vphi=\Dual0\mu{j_{G_0}^\sharp(L_{S(u)}\vphi)}=(-1)^{\Abs0u\Abs0\vphi}\Dual1{\mu_g}{\vphi(\Ad(g^{-1})(u);g)}
  \end{equation}
  \fa open $U\subseteq G_0$, $u\in\Uenv0{\ger g}$, $\mu\in\sh Db_{G_0}(U)$, and $\vphi\in\Gamma_c(\sh O_G|_U)$. 
\end{Prop}

\begin{proof}
  First, we check that the map is well-defined. Indeed, we compute for $x\in\ger g_\ev$:
  \begin{align*}
    \Dual0{ux\otimes\mu}\vphi&=\Dual1\mu{j^\sharp(L_{S(ux)}\vphi)}=-\Dual1\mu{j^\sharp(L_xL_{S(u)}\vphi)}\\
    &=\Dual1{L_x\mu}{j^\sharp{L_{S(u)}\vphi}}=\Dual0{u\otimes L_x\mu}\vphi,
  \end{align*}
  where we abbreviate $j=j_{G_0}$. Similarly, one verifies that the map is $G$-equivariant. Since it is manifestly right $\sh O_{G_0}$-linear, it is a morphism of $(G,\sh O_{G_0})$-modules.

  To see that it is an isomorphism, we define an involutive anti-automorphism $(-)^\vee=i^\sharp:i_0^{-1}\sh O_G\longrightarrow\sh O_G$ (where $i_0(g)=g^{-1}$) by
  \[
    f^\vee(u;g)\defi f\Parens1{\Ad(g)(S(u));g^{-1}}.
  \]
  (This just the inversion morphism $i:G\longrightarrow G$.) Then we compute
  \begin{equation}\label{eq:check-eqn}
    (-1)^{\Abs0u\Abs0\vphi}\Dual0{S(u)\otimes\mu}{\check\vphi}=\Dual1{\mu_g}{\check\vphi\Parens1{\Ad(g^{-1})(u);g}}=\Dual1{\mu_{g^{-1}}}{\vphi(u;g)}.    
  \end{equation}
  We recall again that there is an isomorphism of right $\ger g_\ev$-modules $\Uenv0{\ger g}\cong\textstyle\bigwedge\ger g_\odd\otimes\Uenv0{\ger g_\ev}$. One choice of such an isomorphism that we will repeatedly use is based on the supersymmetrisation map 
  \[
    \beta:S(\ger g)\longrightarrow\Uenv0{\ger g}.
  \]
  Explicitly, it is given by 
  \[
    \textstyle\bigwedge\ger g_\odd\otimes\Uenv0{\ger g_\ev}\longrightarrow\Uenv0{\ger g}:\eta\otimes u\longmapsto\beta(\eta)u,
  \]
  compare \cite{koszul}*{Lemma 1}. Applying this decomposition in Equation \eqref{eq:check-eqn} readily implies our claim. 
\end{proof}

The proof shows that 
\[
  \Uenv0{\ger g}\otimes_{\Uenv0{\ger g_\ev}}\sh Db_{G_0}\cong\textstyle\bigwedge\ger g_\odd\otimes_\knums\sh Db_{G_0}.
\]
If we consider on this sheaf the obvious tensor product locally convex topology (there is no choice which one to take, since $\bigwedge\ger g_\odd$ is finite-dimensional), then it is easy to check that the above isomorphism is in fact one of sheaves of locally convex super-vector spaces. 

\subsection{Left-invariant superdistributions}

In this section, we show that left-in\-va\-ri\-ant superdistributions are smooth and hence proportional to the invariant Berezinian density. To state this precisely, we recall the definition of Berezinian densities. 

\begin{Def}[ber-dens][Berezinian densities]
  Let $\sh Ber_G$ denote the Berezian sheaf of $G$, compare \citelist{\cite{leites}*{2.4.2} \cite{manin}*{4.3.7, 4.6.1} \cite{deligne-morgan}*{\S~1.11, \S~3.10}} for the definition. We let $\Abs0{\Omega}_G\defi or_{G_0}\otimes_\ints\sh Ber_G$, where $or_{G_0}$ is the orientation sheaf of $G_0$. The set of global sections of $\Abs0{\sh Ber}_G$ is denoted by $\Abs0{\Omega}(G)$; elements thereof are called \Define{Berezinian densities}. The set of compactly supported sections of $\Abs0{\sh Ber}_G$ is denoted by $\Abs0{\Omega}_c(G)$.

  Then $\Abs0{\sh Ber}_G$ is naturally a $(G,\sh O_G)$-module. Moreover, if $U\subseteq G_0$ is open and $\omega\in\Gamma_c(\Abs0{\Omega}_G|_U)$, then $\int_{G|_U}\omega\in\knums$, the \Define{Berezin integral} of $\omega$, is well-defined \citelist{\cite{leites}*{Theorem 2.4.5} \cite{manin}*{Theorem 4.6.3} \cite{deligne-morgan}*{Proposition 3.10.5}}. In particular, there is an embedding $\Abs0{\Omega}_G\longrightarrow\sh Db_G$, given by 
  \[
    \Dual0\omega\vphi\defi\int_{G|_U}\omega\vphi\mathfa\omega\in\Abs0{\Omega}_G(U),\vphi\in\Gamma_c(\sh O_G|_U).
  \]
  By \cite{ah-berezin}*{Theorem 4.13}, $\Abs0{\Omega}_G$ has a nowhere vanishing $G$-invariant section $\Abs0{Dg}$, which is unique up to constant multiples. It furnishes a module basis of $\Abs0{\Omega}_G$.
\end{Def}

The following generalises a result due to Bruhat \cite{bruhat}*{Chapitre I, Proposition 3.1}.

\begin{Prop}[bruhat][super Bruhat regularity theorem]
  Let $\mu\in\sh D'(G)$ be left-invariant under $G$. Then for some constant $c$, we have $\mu=c\Abs0{Dg}$.
\end{Prop}

The \emph{proof} of the proposition uses the following definition and basic lemmas.

\begin{Def}[conv-def][convolution of superdistributions]
  Let $\mu,\nu\in\Gamma(\sh Db_G)$. We say that $(\mu,\nu)$ is a \Define{proper pair} if $m_0:\supp\mu\times\supp\nu\longrightarrow G_0$ is a proper map. 

  If $\vphi\in\Gamma_c(\sh O_G)$, then $K\defi(\supp\mu\times\supp\nu)\cap m_0^{-1}(\supp\vphi)$ is compact. Let $\chi\in\Gamma_c(\sh O_{G\times G})$ \scth $\chi|_U=1$ \fs open neighbourhood of $K$. The quantity 
  \[
    \Dual1{\mu*\nu}\vphi\defi\Dual0{\mu\otimes\nu}{\chi m^\sharp(\vphi)}
  \]
  is independent of $\chi$. Moreover, it depends continuously on $\vphi$, thus defining an element $\mu*\nu\in\Gamma(\sh O_G)$, the \Define{convolution} of $\mu$ and $\nu$. Clearly, if either $\mu$ or $\nu$ is compactly supported, then $(\mu,\nu)$ is a proper pair. 
\end{Def}

\begin{Lem}[dist-dens-conv]
    Let $\mu\in\Gamma(\sh D_G)$ and $\omega\in\Abs0\Omega_c(G)$. Then $\mu*\omega\in\Abs0\Omega(G)$.
\end{Lem}

\begin{proof}
  Let $\vphi\in\Gamma_c(\sh O_G)$ and set 
  \[
    \psi(g)\defi\int_G\omega(h)m^\sharp(\vphi)(g,h)
  \]
  for any $T$ and any $h\in_TG$. In the integral, $h$ denotes the generic point $h=\id_G\in_GG$. 

  Then by Yoneda's lemma, we have $\psi\in\Gamma(\sh O_G)$, and this superfunction has compact support $\subseteq(\supp\vphi)(\supp\omega)^{-1}$. Hence, we find that 
  \[
    \Dual1{\mu*\omega}{\vphi}=\Dual1\mu{\psi}=\Dual3{\mu_g}{\int_G\Abs0{Dh}\,f(h)\vphi(gh)}.
  \]

  Writing $\omega=\Abs0{Dg}\,f$, we have 
  \[
    \vrho\defi (\mu\otimes{\id})(m\circ(i\times{\id}))^\sharp(f)\in\Gamma(\sh O_G),
  \]
  since $\Gamma(\sh O_{G\times G})=\Gamma(\sh O_G)\mathop{\widehat{\otimes}}_\pi\Gamma(\sh O_G)$, the completion of the projective tensor product \cite{as}*{Corollary C.9}, and $(\mu\otimes{\id})$ extends continuously this space. We thus compute 
  \[
    \Dual1{\mu*\omega}{\vphi}=\Dual3{\mu_g}{\int_G\Abs0{Dh}\,f(g^{-1}h)\vphi(h)}=\Dual1{\Abs0{Dg}\,\vrho}\vphi,       
  \]
  so that $\mu*\omega=\Abs0{Dg}\,\vrho$, proving the claim. 
\end{proof}

\begin{Lem}[approxid]
  Let $\sh U$ be the filter of open neighbourhoods of $1\in G_0$. There exist Berezinian densities $\chi_U=\check\chi_U\in\Abs0\Omega_c(G)$, $\supp\chi_U\subseteq U\in\sh U$, \scth 
  \[
    \lim_{U\in\mathcal U}\chi_U*\mu=\lim_{U\in\mathcal U}\mu*\chi_U=\mu
  \]
  in $\sh D'(G)$, for any $\mu\in D'(G)$. If $\mu\in\Abs0\Omega_c(G)$, then the convergence is in $\Abs0\Omega_c(G)$.
\end{Lem}

\begin{proof}
  For $U$ sufficiently small, we may choose local coordinates $(u,\xi)$ and define 
  \[
    \chi_U\defi\Abs0{D(u,\xi)}\,\xi_1\dotsm\xi_q\vrho_U,
  \]
  where $\int_U\Abs0{du_0}\vrho_U=1$ and $\dim G=p|q$. Then for $\vphi\in\Gamma_c(\sh O_G)$, we have 
  \[
    \int_G\chi_U\vphi=\int_U\Abs0{du_0}\,\vrho_U\vphi_0\longrightarrow\vphi_0(1)=\vphi(1),
  \]
  where the convergence is uniform for $\vphi$ in compact subsets of $\Gamma_c(\sh O_G)$. Indeed, \cite{fol-abs}*{Proposition 2.42} gives uniform convergence, and compactness is preserved when passing to a coarser topology.

  Now, the computation in the proof of \thmref{Lem}{dist-dens-conv} shows that 
  \[
    \Dual1{\mu*\chi_U}\vphi=\Dual1\mu{\chi_U*\vphi}
  \]
  where we set
  \[
    (\chi_U*\vphi)(h)\defi\int_G\chi_U(g)\vphi(g^{-1}h)
  \]
  for any $T$ and any $h\in_TG$. Then for $h\in_TG$
  \[
    (\chi_U*\vphi)(h)-\vphi(h)=\int_G\chi_U(g)\Parens1{\vphi(g^{-1}h)-\vphi(h)}\longrightarrow0,
  \]
  the convergence being in $\Gamma_c(\sh O_T)$. Taking $T=G$ and $h=\id_G\in_GG$, the assertion follows for right convolutions, and the case of left convolutions is similar. 
\end{proof}

\begin{proof}[\prfof{Prop}{bruhat}]
  The proof is the same as Bruhat's, based on the super-extensions of classical facts stated as the lemmas above. Let $(\chi_U)$ be as in the statement of \thmref{Lem}{approxid}. We have, for any $U\in\sh U$, $g\in G_0$, and $u\in\Uenv0{\ger g}$:
  \[
    L_gL_u(\mu*\chi_U)=(L_gL_u\mu)*\chi_U, 
  \]
  so the superdistribution $\mu*\chi_U$, which is a Berezin density by \thmref{Lem}{dist-dens-conv}, is left-invariant under $G$ and thus equals $c_U\Abs0{Dg}$ for some constant $c_U$. 

  But by \thmref{Lem}{approxid}, we have $\mu=\lim_U\mu*\chi_U=\lim_Uc_U\Abs0{Dg}$, so that $\mu$ is contained in the closure of the line spanned by $\Abs0{Dg}$. But this line is finite-dimensional, and hence a closed subspace of $\Gamma(\sh Db_G)$, since the unique Hausdorff vector space topology on $\knums$ is complete. This shows the assertion. 
\end{proof}

\subsection{Berezinian densities \via ordinary densities}

In this subsection, we show how Berezinian densities can be expressed in terms of ordinary densities on the underlying Lie groups. 

To that end, let $\Abs0{\Omega}_{G_0}\defi or_{G_0}\otimes_\ints\Omega^p_{G_0}$, where $p|q=\dim G$, denote the sheaf of $\knums$-valued smooth densities on $G_0$. Its global sections will be denoted by $\Abs0{\Omega}(G_0)$, and the subspace of compactly supported sections by $\Abs0{\Omega}_c(G_0)$. As above, there is an embedding $\Abs0{\Omega}_{G_0}\longrightarrow\sh Db_{G_0}$, given by 
\[
  \Dual0\omega\vphi\defi\int_U\omega\vphi\mathfa\omega\in\Abs0{\Omega}_{G_0}(U),\vphi\in\Gamma_c(\sh O_{G_0}|_U)=\Ct[^\infty_c]0U.
\]

The isomorphism in \thmref{Prop}{superdist} suggests that we can indentify $\Abs0{\Omega}_G$ and $\Uenv0{\ger g}\otimes_{\Uenv0{\ger g_\ev}}\Abs0{\Omega}_{G_0}$ within $\sh Db_G$. Although this is not completely straightforward, it turns out to be quite generally true, as we now proceed to explain. 

\medskip\noindent
Let $\delta_\odd$ be the character by which $\ger g_\ev$ acts on $\Ber0{\ger g/\ger g_\ev}$, \ie 
\[
  \delta_\odd(x)=-\tr_{\ger g_\odd}\ad(x)\mathfa x\in\ger g_\ev.
\]
This character extends naturally to $\Uenv0{\ger g_\ev}$. It is the differential of the character $\Delta_\odd$ of $G_0$, given by 
\[
  \Delta_\odd(g)\defi\Ber[_{(\ger g/\ger g_\ev)^*}]0{\Ad(g)}=\Parens1{\det\nolimits_{\ger g_\odd}\Ad(g)}^{-1}.
\]

For any $\ger g_\ev$-module $N$ (say), there is a well-known \cites{bell-farnsteiner,chemla-gro,chemla-dual,dp,gorelik-ghost} isomorphism of graded $\ger g$-modules
\begin{center}
  \begin{tikzcd}
    \Uenv0{\ger g}\otimes_{\Uenv0{\ger g_\ev}}N\rar{\Phi}&
    \GHom[_{\ger g_\ev}]1{\Uenv0{\ger g},\Ber0{(\ger g/\ger g_\ev)^*}\otimes_\knums N}.
  \end{tikzcd}
\end{center}
By the construction detailed in \cite{gorelik-ghost}*{3.2.1}, it is given explicitly by 
\begin{equation}\label{eq:phidef}
  \Phi(u\otimes n)(v)=(-1)^{(\Abs0u+\Abs0n)\Abs0v}\iota(vu)(\omega_\odd\otimes n)
\end{equation}
where $\omega_\odd\in\Ber0{(\ger g/\ger g_\ev)^*}$ is an arbitrary non-zero element and $\iota:\Uenv0{\ger g}\longrightarrow\Uenv0{\ger g_\ev}$ is the left $\ger g_\ev$-linear map defined by
\[
  \iota(u\beta(\eta))\defi u\int_{\ger g_\odd}\omega_\odd\eta\mathfa u\in\Uenv0{\ger g_\ev},\eta\in\textstyle\bigwedge\ger g_\odd.
\]
Here, $\beta$ is supersymmetrisation, and the Berezin integral is normalised by $\int_{\ger g_\odd}\omega_\odd=1$. 

A notable special case occurs when $N=\Ber0{\ger g/\ger g_\ev}$. In this case, we may consider the action of $G_0$ on $N$, and $\Ber0{\ger g/\ger g_\ev}^*\otimes_\knums N\cong\knums$ as $G$-modules. Moreover, if $g\in G_0$, then we have the equation 
\[
  \iota\Parens1{\Ad(g)(u\beta(\eta))}=\Ad(g)(u)\int_{\ger g_\odd}\omega_\odd\Ad(g)(\eta)=\Delta_\odd(g)\cdot\Ad(g)\Parens1{\iota(u\beta(\eta))}
\]
by the change of variables formula for the Berezin integral. Combining these facts with the definition of $\Phi$, one arrives by the formula
\begin{equation}\label{eq:phi-ad}
  \Phi\Parens1{\Ad(g)(u)\otimes n}(v)=\Delta_\odd(g)\Phi(u\otimes n)\Parens1{\Ad(g^{-1})(v)}
\end{equation}
for $u,v\in\Uenv0{\ger g}$ and $n\in\Ber0{\ger g/\ger g_\ev}$.

Let $I_{\delta_\odd}$ be the left ideal of $\Uenv0{\ger g}$ generated by the set
\[
  \Set1{x\in\ger g_\ev}{x-\delta_\odd(x)}.
\]
By \cite{dp}*{Proposition 3.5}, the space of $\ger g$-invariants in 
\[
  \Uenv0{\ger g}\otimes_{\Uenv0{\ger g_\ev}}\Ber0{\ger g/\ger g_\ev}=\Uenv0{\ger g}/I_{\delta_\odd}
\]
is one-dimensional. Let $\gamma\in\Uenv0{\ger g}$ be a representative of a basis. 

\begin{Prop}[invber]
  For a suitable normalisation of $\Abs0{Dg}$ and $\Abs0{dg}$, we have
  \begin{equation}\label{eq:invber}
    \Abs0{Dg}=L_\gamma(\Abs0{dg}\,\Delta_\odd).
  \end{equation}
\end{Prop}

\begin{proof}
  Let us consider the isomorphism $\Phi$ for $N=\Ber0{\ger g/\ger g_\ev}$. Since there is a canonical isomorphism
  \begin{equation}\label{eq:bertensor}
    \Ber0{(\ger g/\ger g_\ev)^*}\otimes_\knums\Ber0{\ger g/\ger g_\ev}\longrightarrow\knums
  \end{equation}
  of $\ger g_\ev$-modules \cite{dp}*{Lemma 1.4}, we may view $\Phi$ as an isomorphism
  \[
    \Phi:\Uenv0{\ger g}\otimes_{\Uenv0{\ger g_\ev}}\Ber0{\ger g/\ger g_\ev}\longrightarrow\GHom[_{\ger g_\ev}]0{\Uenv0{\ger g},\knums}.
  \]

  Moreover, by \cite{dp}*{p.~150}, the coset of $\gamma$ corresponds under the canonical isomorphism $\Phi$ to the element $\eps:\Uenv0{\ger g}\longrightarrow\knums$, which is the extension of $0:\ger g\longrightarrow\knums$ to a superalgebra morphism. Hence, by Equation \eqref{eq:phi-ad}, for any $g\in G_0$, we have 
  \begin{equation}\label{eq:ad-gamma}
    \Ad(g)(\gamma)\equiv\Delta_\odd(g)\gamma\pmod{I_{\delta_\odd}}.
  \end{equation}
  
  Since $\Abs0{dg}\,\Delta_\odd$ is relatively $\ger g_\ev$-invariant for the character $\delta_\odd$, this quantity is annihilated by $I_{\delta_\odd}$. In particular, the superdistribution 
  \[
    \Omega\defi L_\gamma(\Abs0{dg}\,\Delta_\odd)\in\sh D'(G)
  \]
  depends only on the coset of $\gamma$. 

  By \thmref{Prop}{bruhat}, it will be sufficient to show that $\Omega$ is a $\ger g$- and $G_0$-invariant functional. First, let $x\in\ger g$ be homogeneous. Then we compute 
  \[ 
    \Dual0\Omega{L_xf}=\Dual1{\Abs0{dg}\,\Delta_\odd}{L_{S(\gamma)x}(f)}
    =(-1)^{\Abs0x\Abs0\gamma}\Dual1{L_{x\gamma}(\Abs0{dg}\,\Delta_{\delta_\odd})}{f}=0
  \]
  since by the choice of $\gamma$, we have $x\gamma\in I_{\delta_\odd}$ for any $x\in\ger g$. 

  Secondly, we compute 
  \[
    \Dual0\Omega{L_hf}=\Dual1{L_{\Ad(h^{-1})(\gamma)}(L_{h^{-1}}(\Abs0{dg}\,\Delta_\odd))}{f}=\Dual0\Omega f,
  \]
  by the use of the relation $L_{h^{-1}}(\Abs0{dg}\,\Delta_\odd)=\Delta_\odd(h)\,\Abs0{dg}\,\Delta_\odd$ and Equation \eqref{eq:ad-gamma}. Thus, we reach our conclusion.
\end{proof}

\begin{Cor}[ber-dens]
  As $(G,\sh O_{G_0})$-submodules of $\sh Db_G$, we have 
  \[
    \Abs0{\Omega}_G=\Uenv0{\ger g}\otimes_{\Uenv0{\ger g_\ev}}\Abs0{\Omega}_{G_0}. 
  \]
  Writing $\Delta(\gamma)=\sum_i\gamma_i'\otimes\gamma_i''$, the Berezinian density $\Abs0{Dg}\,f$ corresponds to 
  \begin{equation}\label{eq:invber-haar}
    \sum\nolimits_i\gamma_i'\otimes\Abs0{dg}\,\Delta_\odd\,j_{G_0}^\sharp\Parens1{L_{S(\gamma_i'')}(f)}.
  \end{equation}
  Conversely, the element $1\otimes\Abs0{dg}$ is mapped to $\Abs0{Dg}\,\psi$, where $\psi\in\Gamma(\sh O_G)$ is defined by 
  \begin{equation}\label{eq:psidef}
    \psi(u;g)\defi (R_{\iota(u)}\Delta_\odd^{-1})(g)
  \end{equation}
  \fa $u\in\Uenv0{\ger g}$, $g\in G_0$. 
\end{Cor}

\begin{proof}
  Consider the isomorphism 
  \begin{center}
    \begin{tikzcd}
      \Uenv0{\ger g}\otimes_{\Uenv0{\ger g_\ev}}\sh Db_{G_0}\rar{}&\sh Db_G
    \end{tikzcd}
  \end{center}
  from \thmref{Prop}{superdist}. For $f,\vphi\in\sh O_G(U)$, we compute
  \[
    \Parens1{L_{S(\gamma)}(f\vphi)}(g)=\sum\nolimits_i(-1)^{\Abs0{\gamma_i'}(\Abs0f+\Abs0{\gamma_i''})}\Parens1{L_{S(\gamma_i'')}(f)}(g)\Parens1{L_{S(\gamma_i')}(\vphi)}(g).
  \]
  For the non-zero summands, we have $\Abs0{\gamma_i''}+\Abs0f\equiv0\,(2)$. Hence, under the isomorphism, the expression in Equation \eqref{eq:invber-haar} is mapped to $\Abs0{Dg}\,f$. Thus, $\Abs0{\Omega}_G$ is contained in the image of the subsheaf
  \[
    \Uenv0{\ger g}\otimes_{\Uenv0{\ger g_\ev}}\Abs0{\Omega}_{G_0}\subseteq\Uenv0{\ger g}\otimes_{\Uenv0{\ger g_\ev}}\sh Db_{G_0}.
  \]

  For the converse, \ie that $\Uenv0{\ger g}\otimes_{\ger g_\ev}\Abs0{\Omega}_{G_0}$ is mapped to $\Abs0{\Omega}_G$, we need only show that this is the case for the $\Uenv0{\ger g}\otimes\sh O_{G_0}$-generator $1\otimes\Abs0{dg}$. To that end, consider the superfunction $\psi\in\Gamma(\sh O_G)$, defined by Equation \eqref{eq:psidef}. It is well-defined, because the map $\iota$ is by definition left $\ger g_\ev$-linear.
  
  By \cite{dp}*{Theorem 3.1, Equation (65)}, we have $\gamma\equiv\beta(x_1\dotsm x_qJ)\ (I_{\delta_\odd})$, where $J\in(\bigwedge\ger g_\odd)_\ev$ is the Jacobian of the exponential map (compare \loccit). Set $\tilde\gamma\defi x_1\dotsm x_qJ\in S(\ger g)$ and consider the grading with components 
  \[
    S^{\bullet,k}\defi S^{\bullet,k}(\ger g)\defi S(\ger g_\ev)\otimes\textstyle\bigwedge^k\ger g_\odd.  
  \] 
  Observe 
  \[
    \Delta(S^{\bullet,k})\subseteq\textstyle\bigoplus_{a+b=k}S^{\bullet,a}\otimes S^{\bullet,b}.
  \]
  In particular, we have
  \[
    \Delta(\tilde\gamma)\equiv\tilde\gamma\otimes1\pmod{\textstyle\bigoplus_{a<q}S^{\bullet,q}\otimes S(\ger g)}.
  \]
  On the other hand, by the definition of $\iota$, we have $\iota(u\beta(\eta))=0$ for $u\in\Uenv0{\ger g_\ev}$ and $\eta\in\bigwedge\ger g_\odd$, unless $\eta$ has a non-zero component in top degree. Since $\beta:\Uenv0{\ger g}\longrightarrow S(\ger g)$ is an isomorphism of coalgebras \cite{pet}*{Theorem 8.1}, we find
  \[
    \Delta(\gamma)\equiv\gamma\otimes 1\pmod{\ker\iota\otimes\Uenv0{\ger g}}.
  \]

  As observed in the proof of \thmref{Prop}{invber}, we have \fa $v\in\Uenv0{\ger g}$:
  \[
    \eps(v)=\Phi(\dot\gamma)(v)=\iota(v\gamma),
  \]
  where $\dot\gamma\in\Uenv0{\ger g}/I_{\delta_\odd}$ denotes the coset of $\gamma$ and we have used Equation \eqref{eq:phidef}. No signs occur, since the left-hand side of the equation is independent of the odd part of $v$. In particular, $\iota(\gamma)=1$. 

  Hence, we compute for any $\vphi\in\Gamma(\sh O_G)$ that 
  \[
    L_{S(\gamma)}(\psi\vphi)(1;g)=\Delta_\odd(g)^{-1}\vphi(1;g).
  \]
  For compactly supported $\vphi$, this implies that 
  \begin{align*}
    \Dual1{\Abs0{Dg}\,\psi}\vphi=\Dual1{\Abs0{dg}}{\Delta_\odd j_{G_0}^\sharp\Parens1{L_{S(\gamma)}(\psi\vphi)}}=\Dual1{\Abs0{dg}}{j_{G_0}^\sharp(\vphi)}.
  \end{align*}
  Thus, we find that $1\otimes\Abs0{dg}$ is mapped to $\Abs0{\Omega}(G)$; this proves the claim.
\end{proof}

\subsection{Convolution of superdistributions and Berezinian densities}

\begin{Def}[conv-dist][compactly supported superdistributions]
  We let $\sh E'(G)$ be the strong dual space of $\Gamma(\sh O_G)=\sh O_G(G_0)$ and call its elements \Define{compactly supported superdistributions}. For $\mu,\nu\in\sh E'(G)$, the convolution $\mu*\nu\in\sh E'(G)$ from \thmref{Def}{conv-def} takes the form 
  \[
    \Dual0{\mu*\nu}f\defi\Dual0{\mu\otimes\nu}{m^\sharp(f)}
  \]
  \fa $f\in\Gamma(\sh O_G)$. Here, $m:G\times G\longrightarrow G$ is the multiplication of $G$.

  If $A$ is a topological $\knums$-vector space with an algebra structure, then we call $A$ a \Define{topological algebra} if multiplication is separately continuous. We allow non-unital algebras, but unless called `non-unital' expressly, they are assumed to have a unit.
\end{Def}

In the following, let $\sh E'(G_0)$ be the strong dual of $\Gamma(\sh O_{G_0})$. It carries a natural convolution, see \citelist{\cite{bruhat}*{\S 1.4} \cite{kv-coho}*{Chapter I.1} \cite{taylor-ncharman}*{Chapter 0.3}}. Recall that $\Uenv0{\ger g_\ev}\subseteq\sh E'(G_0)$ is a subalgebra \via $u\longmapsto L_u\delta$, where $\delta$ denotes the Dirac delta distribution supported at the neutral element of $G_0$.

\begin{Prop}[dist-conv]
  The convolution product on $\sh E'(G)$ is well-defined and turns it into an associative and unital topological superalgebra. We have $\sh E'(G)=\Gamma_c(\sh Db_G)$ and there is an isomorphism 
  \begin{center}
    \begin{tikzcd}
      \Uenv0{\ger g}\otimes_{\Uenv0{\ger g_\ev}}\sh E'(G_0)\rar{}&\sh E'(G)
    \end{tikzcd}
  \end{center}
  of locally convex super-vector spaces. In terms of this isomorphism, the superalgebra structure is uniquely determined by the following facts:
  \begin{enumerate}[wide]
    \item\label{dist-conv-i} The following are graded subalgebras:
    \[
      \Uenv0{\ger g}=\Uenv0{\ger g}\otimes_{\Uenv0{\ger g_\ev}}\Uenv0{\ger g_\ev},\quad
      \sh E'(G_0)=\Uenv0{\ger g_\ev}\otimes_{\Uenv0{\ger g_\ev}}\sh E'(G_0).
    \]
    \item\label{dist-conv-ii} For $\mu\in\sh E'(G_0)$ and $u\in\Uenv0{\ger g}$, the products $u*\mu$ and $\mu*u$ are given by 
    \begin{equation}\label{eq:conv}
      \Dual1{u*\mu}\vphi=\Dual1\mu{j_{G_0}^\sharp\Parens1{L_{S(u)}\vphi}},\quad
      \Dual1{\mu*u}\vphi=\Dual1\mu{j_{G_0}^\sharp\Parens1{R_u\vphi}}
    \end{equation}
    \fa superfunctions $\vphi\in\Gamma(\sh O_G)$.
  \end{enumerate}
\end{Prop}

\begin{proof}
  Since $\Gamma_c(\sh O_G)$ is dense in $\Gamma(\sh O_G)$, $\sh E'(G)$ may be identified with a subspace of $\sh D'(G)$. On the other hand, one knows that $\sh E'(G_0)=\Gamma_c(\sh Db_{G_0})$. Therefore, \thmref{Prop}{superdist} gives an isomorphism of super-vector spaces as stated and $\sh E'(G)=\Gamma_c(\sh Db_G)$. Moreover, it is straightforward to prove that it is indeed a homeomorphism for the topology on $\sh E'(G)$ and the natural topology on $\Uenv0{\ger g}\otimes_{\Uenv0{\ger g_\ev}}\sh E'(G_0)=\bigwedge\ger g_\odd\otimes_\knums\sh E'(G_0)$. 

  It is clear that there is at most one algebra structure on $\sh E'(G)$ determined by the information stated in \eqref{dist-conv-i} and \eqref{dist-conv-ii}. Conversely, we compute for $\mu,\nu\in\sh E'(G_0)$ and $u,v\in\Uenv0{\ger g}$, by the use of Equations \eqref{eq:mult} and \eqref{eq:tensor-dist}:
  \begin{align*}
    \Dual1{(u\otimes\mu)*(v\otimes\nu)}{\vphi}&=(-1)^{\Abs0\vphi(\Abs0u+\Abs0v)}\Dual1{\mu_g\otimes\nu_h}{(m^\sharp\vphi)\Parens1{\Ad((g,h)^{-1})(u\otimes v);g,h}}\\
    &=(-1)^{\Abs0\vphi(\Abs0u+\Abs0v)}\Dual1{\mu_g\otimes\nu_h}{\vphi\Parens1{\Ad(h^{-1})(\Ad(g^{-1})(u)v);gh}}.
  \end{align*}
  For $\mu=\delta$ and $v=1$, we obtain 
  \[
    \Dual1{u*\nu}\vphi=(-1)^{\Abs0\vphi\Abs0u}\Dual1{\nu_h}{\vphi\Parens1{\Ad(h^{-1})(u);h}}=\Dual1\nu{L_{S(u)}\vphi},
  \]
  and for $u=1$ and $\nu=\delta$, we get  
  \[
    \Dual1{\mu*v}\vphi=(-1)^{\Abs0\vphi\Abs0v}\Dual1{\mu_g}{\vphi\Parens0{v;g}}=\Dual1\mu{R_v\vphi}.
  \]
  This shows Equation \eqref{eq:conv}. 

  The convolution on $\sh E'(G)$ is an even bilinear map by definition. That it is an associative operation follows either from $m\circ(m\times\id)=m\circ(\id\times m)$, or also easily from Equation \eqref{eq:conv}, together with the fact that $\Uenv0{\ger g}$ and $\sh E'(G_0)$ are algebras and that the actions $L$ and $R$ commute. 
\end{proof}

The convolution algebra structure on $\sh E'(G)$ admits a natural $\knums$-linear anti-involution, defined by 
\begin{equation}\label{eq:checkdef}
  \Dual0{\check\mu}{\vphi}\defi\Dual0\mu{\check\vphi}=\Dual0\mu{i^\sharp\vphi},  
\end{equation}
where $i:G\longrightarrow G$ is the inversion morphism, and $\check\vphi=i^\sharp\vphi$ was employed above in the proof of \thmref{Prop}{superdist}. Since Berezinian densities pull back under isomorphisms, the involution leaves $\Abs0\Omega_c(G)\subseteq\sh E'(G)$ stable. 

\begin{Cor}[dens-conv]
  The dense subspace $\Abs0\Omega_c(G)\subseteq\sh E'(G)$ is a graded ideal and a non-unital Fr\'echet algebra with the topology induced from $\Gamma_c(\sh O_G)$. In terms of the isomorphism $\Abs0\Omega_c(G)=\Uenv0{\ger g}\otimes_{\Uenv0{\ger g_\ev}}\Abs0\Omega_c(G_0)$, its $\sh E'(G)$-bimodule structure is determined uniquely by the following facts:
  \begin{enumerate}[wide]
    \item\label{dens-conv-i} The following is a non-unital graded subalgebra bi-invariant under $\sh E'(G_0)$:
    \[
      \Abs0\Omega_c(G_0)=\Uenv0{\ger g_\ev}\otimes_{\Uenv0{\ger g_\ev}}\Abs0\Omega_c(G_0).
    \]
    \item\label{dens-conv-ii} For $u,v\in\Uenv0{\ger g}$ and $\omega\in\Abs0\Omega_c(G_0)$, we have
    \[
      u*(v\otimes\omega)=(u\otimes1)*(v\otimes\omega)=uv\otimes\omega.
    \]
    \item\label{dens-conv-iii} For $\omega\in\Abs0\Omega_c(G_0)$ and $u\in\Uenv0{\ger g}$, the products $u*\omega$ and $\omega*u$ are given by 
    \begin{equation}\label{eq:dens-conv}
      \begin{aligned}[c]
        \int_{G_0}(u*\omega)\vphi&=\int_{G_0}\omega\,j_{G_0}^\sharp\Parens1{L_{S(u)}\vphi}=\Dual1{u\otimes\omega}\vphi,\\
        \int_{G_0}(\omega*u)\vphi&=\int_{G_0}\omega\,j_{G_0}^\sharp\Parens1{R_u\vphi},
      \end{aligned}
    \end{equation}
    \fa superfunctions $\vphi\in\Gamma(\sh O_G)$.
  \end{enumerate}
\end{Cor}

\begin{proof}
  Let us verify that $\Abs0\Omega_c(G)$ is indeed a convolution ideal in $\sh E'(G)$. Indeed, this follows from \thmref{Lem}{dist-dens-conv}. Alternatively, one may proceed as follows. 

  Certainly, $\Abs0\Omega_c(G_0)$ is an ideal of $\sh E'(G_0)$. Let $\omega$ and $u\in\Uenv0{\ger g}$. Since $u*\omega$ corresponds to $u\otimes\omega$, it is obvious that $u*\omega\in\Omega_c(G)$. On the other hand, we have 
  \[
    \int_{G_0}(\omega*u)\vphi=(-1)^{\Abs0u\Abs0\vphi}\int_{G_0}\omega(g)\vphi(u;g)=(-1)^{\Abs0u\Abs0\vphi}\int_{G_0}\check\omega(g)\check\vphi(S(u);g)=\int_G\check\Omega\,\vphi,
  \]
  where $\Omega\in\Abs0\Omega_c(G)$ corresponds to $S(u)\otimes\check\omega$ and $\check\Omega$ was defined in Equation \eqref{eq:checkdef}. This shows that $\omega*u\in\Abs0\Omega_c(G)$.

  Thus, $\Abs0\Omega_c(G)$ is indeed a graded ideal of $\sh E'(G)$, and the remaining statements follow readily from \thmref{Prop}{dist-conv}.
\end{proof}

\subsection{Convolution action on representations}

We now show how supergroup representations on Fr\'echet spaces can be characterised in terms of the action of convolution superalgebras. We will use the following terminology.

\begin{Def}
   Left $A$ be a topological algebra. A left $A$-module will be called a \Define{continuous module} if the action map is separately continuous. An $A$-module $E$ is called \Define{non-degenerate} if 
  \[
    E=AE\defi\Span1{av\bigm|a\in A,v\in E}_\knums.
  \]
\end{Def}

\begin{Lem}[intrep]
  Let $(E,\pi)$ be a weakly smooth Fr\'echet $G$-representation. Then the $\sh E'(G_0)$-module structure inherited from $E|_{G_0}$ combines with the $\Uenv0{\ger g}$-action on $E$ to a unique continuous $\sh E'(G)$-module structure on $E$, denoted by $\Pi$. The action map $\sh E'(G)\times E\longrightarrow E$ is hypocontinuous. 
\end{Lem}

\begin{proof}
  We begin with some preliminary considerations. Since $\Gamma(\sh O_{G_0})$ is a nuclear Fr\'echet space \cite{treves}*{Corollary to Theorem 51.4}, we have $\Ct[^\infty]0{G_0,E}=\Gamma(\sh O_{G_0})\mathop{\widehat{\otimes}}_\pi E$ by \cite{treves}*{Theorems 44.1 and 50.1}, where $\widehat{\otimes}_\pi$ denotes the completed projective tensor product. Moreover, $\Gamma(\sh O_{G_0})$ is barreled and reflexive, and $\sh E'(G_0)$ is nuclear and complete in the strong topology \cite{treves}*{Corollary 2 to Theorem 32.2, Corollary 1 to Proposition 33.2, Corollary to Proposition 36.9, Proposition 36.10, Proposition 50.6}. 

  In particular, the abstract Kernels Theorem \cite{treves}*{Proposition 50.5} applies. Thus, if $\Dual0\cdot\cdot$ denotes the pairing of $\sh E'(G_0)$ and $\Gamma(\sh O_{G_0})$, then the map
  \[
    \Ct[^\infty]0{G_0,E}=\Gamma(\sh O_{G_0})\mathop{\widehat{\otimes}}\nolimits_\pi E=\sh E'(G_0)'_\beta\mathop{\widehat{\otimes}}\nolimits_\pi E\longrightarrow\GHom0{\sh E'(G_0),E}
  \]
  induced by $\vphi\otimes e\longmapsto(\mu\longmapsto\Dual0\mu\vphi e)$ is a continuous linear isomorphism. Here, $\underline{\mathrm{Hom}}$ denotes the space of continuous linear maps with the topology of uniform convergence on bounded subsets. This map is an element of 
  \[
    \GHom1{\Ct[^\infty]0{G_0,E},\GHom0{\sh E'(G_0),E}},
  \]
  so by \cite{bourbaki-evt}*{Chapitre III, \S~5.3, Proposition 3}, the corresponding bilinear map
  \[
    \Dual0\cdot\cdot:\sh E'(G_0)\times\Ct[^\infty]0{G_0,E}\longrightarrow E,
  \]
  is hypocontinuous with respect to the first argument. Since $\sh E'(G_0)$ is barreled, it is hypocontinuous \cite{bourbaki-evt}*{Chapitre III, \S~5.3, Proposition 6}.

  In other words, the bilinear map
  \[
    \sh E'(G_0)\times(\Gamma(\sh O_{G_0})\otimes E)\longrightarrow E
  \]
  that sends $\Parens0{\mu,\vphi\otimes e}$ to $\Dual0\mu\vphi e$ has a (unique) hypocontinuous bilinear extension.
  
  Next, recall that $\Gamma(\sh O_G)\cong\GHom[_{\ger g_\ev}]0{\Uenv0{\ger g},\Gamma(\sh O_{G_0})}$, see the beginning of Subsection \ref{subs:superdist}. Here, $\Gamma(\sh O_{G_0})\otimes\bigwedge(\ger g_\odd)^*$ induces the locally convex topology on the latter space. Thus, we define analogously:
  \[
    \Ct[^\infty]0{G,E}\defi\Gamma(\sh O_G)\mathop{\widehat{\otimes}}\nolimits_\pi E.
  \]
  The above arguments apply, and we get a natural hypocontinuous bilinear map
  \[
    \Dual0\cdot\cdot:\sh E'(G)\otimes\Ct[^\infty]0{G,E}\longrightarrow E.
  \]
  One sees that the maps thus constructed are compatible with the natural maps $\sh E'(G_0)\longrightarrow\sh E'(G)$ and $\Ct[^\infty]0{G,E}\longrightarrow\Ct[^\infty]0{G_0,E}$. We will therefore suppress these in the notation. Moreover, from the definition, it is clear that 
  \begin{equation}\label{eq:ev-ev}
    \Dual1\mu{T\circ\vphi}=T\Parens1{\Dual0\mu\vphi}
  \end{equation}
  for any continuous linear endomorphism $T$ of $E$ and any $\mu\in\sh E'(G)$, $\vphi\in\Ct[^\infty]0{G,E}$. This completes our preliminaries. 

  Now, let $d\pi$ and $\pi_0$ denote the action of $\ger g$ and $G_0$ on $E$, respectively. Take $v\in E$. Then $\pi_0(-)v:G\longrightarrow E:g\longmapsto\pi_0(g)v$ is a smooth map and there is an element $\pi_v=\pi(-)v\in\GHom[_{\ger g_\ev}]0{\Uenv0{\ger g},\Ct[^\infty]0{G_0,E}}$, defined by
  \[
    \pi_v(u;g)\defi\pi(u;g)v\defi\Parens1{\pi(-)v}(u)(g)\defi\pi_0(g)d\pi(u)v.
  \]

  We may thus define for $u\in\Uenv0{\ger g}$ and $\mu\in\sh E'(G_0)$:
  \begin{equation}\label{eq:intrep-def}
    \Pi(u\otimes\mu)v\defi d\pi(u)\Pi_0(\mu)v=d\pi(u)\Dual1{\mu_g}{\pi_0(g)v},
  \end{equation}
  where we let $\Pi_0$ denote the `integrated' version (on distributions) of the $G_0$-representation $\pi_0$ on $E$ \cite{taylor-ncharman}*{Chapter 0.3}. 
  
  We compute with Equation \eqref{eq:ev-ev} and \thmref{Prop}{dist-conv} that 
  \[
    \begin{split}
      \Pi(u\otimes\mu)v&=\Dual1{\mu_g}{d\pi(u)\pi_0(g)v}=\Dual1{\mu_g}{\pi_0(g)d\pi(\Ad(g^{-1})(u))v}\\
      &=(-1)^{\Abs0u\Abs0v}\Dual1\mu{j_{G_0}^\sharp(L_{S(u)}(\pi_v))}=\Dual1{u\otimes\mu}{\pi_v},
    \end{split}
  \]
  \ie $\Pi(\nu)v=\Dual0\nu{\pi_v}$ \fa $\nu\in\sh E'(G)$.

  In particular, for $x\in\ger g_\ev$, we obtain
  \[
    \begin{split}
      \Pi(ux\otimes\mu)v&=-(-1)^{\Abs0u\Abs0v}\Dual1\mu{L_xL_{S(u)}\pi_v}\\
      &=(-1)^{\Abs0u\Abs0v}\Dual1{L_x\mu}{L_{S(u)}\pi_v}=\Pi(u\otimes L_x\mu)v.
    \end{split}
  \]
  This shows that the action $\Pi$ is well-defined on $\Uenv0{\ger g}\otimes_{\Uenv0{\ger g_\ev}}\sh E'(G_0)$, and by our preliminary considerations, the action map $\sh E'(G)\times E\longrightarrow E$ is hypocontinuous.

  To see that $\Pi$ is an action, let $u\in\Uenv0{\ger g}$, $\mu\in\sh E'(G_0)$, and $v\in E$. Then $\Pi(u\otimes1)\Pi(1\otimes\mu)v=\Pi(u\otimes\mu)v$ from the above. Moreover, we have 
  \[
    \pi_{d\pi(u)v}(u';g)=\pi_0(g)d\pi(u')d\pi(u)v=\pi_0(g)d\pi(u'u)v=(-1)^{\Abs0u\Abs0{u'}}(R_u\pi_v)(u';g),
  \]
  so 
  \[
    \begin{split}
      \Pi(1\otimes\mu)\Pi(u\otimes1)v&=\Pi(1\otimes\mu)d\pi(u)v=\Dual1\mu{\pi_{d\pi(u)v}}\\
      &=\Dual1\mu{R_u\pi_v}=\Dual1{\mu*u}{\pi_v}=\Pi(\mu*u)v
    \end{split}
  \]
  by Equation \eqref{eq:conv}. This proves the claim, since $(1\otimes\mu)(u\otimes1)$ and $\mu*u$ are identified within $\sh E'(G)$ by \thmref{Prop}{dist-conv}.
\end{proof}

We call $\Pi$ the \Define{integrated action} of $\pi$. Restricting it to densities, we obtain the following proposition, which generalises a theorem of Dixmier--Malliavin \cite{dix-mal}.

\begin{Prop}[dm][super Dixmier--Malliavin theorem]
  Let $E$ be a Fr\'echet super-vector space over $\knums$. Then we have the following facts: 
  \begin{enumerate}[wide]
    \item\label{dm-i} If $E$ carries the structure of a continuous $G$-representation, then the action of $\Abs0\Omega_c(G)$ on $E_\infty$ extends continuously to $E$. The induced action of $\Abs0\Omega_c(G)$ on $E_\infty$ is non-degenerate. More precisely, we have the equality
    \begin{equation}\label{eq:dm}
      E_\infty=\Pi\Parens1{\Abs0\Omega_c(G)}E=\Pi\Parens1{\Abs0\Omega_c(G)}E_\infty.
    \end{equation}
    \item\label{dm-ii} Conversely, let $\Pi$ be a non-degenerate continuous action of $\Abs0\Omega_c(G)$ on $E$. Then $\Pi$ is integrated from a unique weakly smooth $G$-representation. 
  \end{enumerate}
  In particular, the category of weakly smooth Fr\'echet $G$-representations and the category of non-degenerate continuous Fr\'echet $\Abs0\Omega_c(G)$-modules are equivalent. 
\end{Prop}

\begin{proof}
  Assume that $\pi$ is a continuous $G$-representation on $E$, so that we have by \thmref{Lem}{intrep} the integrated representation $\Pi$ of $\sh E'(G)$ on $E_\infty$. Let $\Pi_0$ be the integrated version of the $G_0$-representation $\pi_0$ of $G_0$ on $E$. By a theorem of G\aa{}rding \cite{taylor-ncharman}*{(2.22)}, we have $\Pi_0\Parens1{\Abs0\Omega_c(G_0)}E\subseteq E_\infty$. 

  Thus, it makes sense to define, for $u\in\Uenv0{\ger g}$, $\omega\in\Abs0\Omega_c(G_0)$, and $v\in E$:
  \[
    \Pi(u\otimes\omega)v\defi d\pi(u)\Pi_0(\omega)v.
  \]
  Indeed, this coincides with the definition of $\Pi$ on $E_\infty$ given in Equation \eqref{eq:intrep-def}. In addition, for $x\in\ger g_\ev$, we have
  \[
    \Pi(ux\otimes\omega)v=d\pi(u)d\pi_0(x)\Pi_0(\omega)v=d\pi(u)\Pi_0(L_x\omega)=\Pi(u\otimes x\omega)v,
  \]
  so that $\Pi$ defines a continuous representation of $\Abs0\Omega_c(G)$ by \thmref{Cor}{ber-dens}.

  By the Dixmier--Malliavin theorem \cite{dix-mal}*{Theorem 3.3}, we have
  \[
    E_\infty=\Pi_0\Parens1{\Abs0\Omega_c(G_0)}E_\infty.
  \] 
  Applying the definition of $\Pi$ in Equation \eqref{eq:intrep-def}, we obtain
  \[
    E_\infty=d\pi(\Uenv0{\ger g})E_\infty=d\pi(\Uenv0{\ger g})\Pi_0\Parens1{\Abs0\Omega_c(G_0)}E_\infty=\Pi\Parens1{\Abs0\Omega_c(G)}E_\infty.
  \]
  But $\Pi(\Abs0\Omega_c(G))E\subseteq E_\infty$, so we have proved part \eqref{dm-i} of the proposition. 

  Conversely, assume that $E$ is a non-degenerate continuous $\Abs0\Omega_c(G)$-module. If $v\in E$ is a vector, then we may express it as $v=\sum_{j\in J}\Pi(\omega_j)v_j$ where $J$ is finite and $\omega_j\in\Abs0\Omega_c(G)$, $v_j\in V$. We wish to define $\pi_0$ and $d\pi$ for $g\in G_0$ and $u\in\Uenv0{\ger g}$ by
  \begin{equation}\label{eq:intrep-to-rep}
    \pi_0(g)v\defi\sum\nolimits_{j\in J}\Pi(L_g\omega_j)v_j,\quad d\pi(u)v\defi\sum\nolimits_{j\in J}\Pi(L_u\omega_j)v_j.
  \end{equation}
  The first task is to show that these quantities are independent of all choices.

  To that end, let $\sum_{j\in J}\Pi(\omega_j)v_j=0$ in $E$. Choose $(\chi_U)$ as in \thmref{Lem}{approxid}. Then $L_g\omega_j=\delta_g*\omega_j$, where $\delta_g$ is the Dirac distribution supported at $g$, and 
  \[
    \sum\nolimits_{j\in J}\Pi(L_g\omega_j)v_j=\lim_{U\in\mathcal U}\sum\nolimits_{j\in J}\Pi(\delta_g*\chi_U*\omega_j)v_j=\lim_{U\in\mathcal U}\Pi(L_g\chi_U)v=0. 
  \]
  A similar argument applies for $u\in\Uenv0{\ger g}$, and so Equation \eqref{eq:intrep-to-rep} indeed defines actions $\pi_0$ of $G_0$ and $d\pi$ of $\ger g$. 

  Moreover, in case $\Pi$ is already integrated from a weakly smooth $G$-representation $\pi'$, then analogously 
  \[
    \pi'_0(g)v=\lim_{U\in\mathcal U}\sum\nolimits_j\pi'_0(g)\Pi(\chi_U*\omega_j)v_j=\lim_{U\in\mathcal U}\Pi(L_g\chi_U)v=\pi_0(g)v.
  \]
  Similarly, one shows that $d\pi'=d\pi$, so $\Pi$ is integrated from at most one weakly smooth $G$-representation, and if it is, then the corresponding actions of $G_0$ and $\ger g$ are given by $\pi_0$ and $d\pi$, respectively. It therefore remains to be shown that $\pi_0$ and $d\pi$ combine to a weakly smooth $G$-representation. 

  For this, we observe that the action map $\Pi:\Abs0\Omega_c(G)\otimes_i E\longrightarrow E$ is continuous and surjective, $\otimes_i$ denoting injective tensor product. It extends to a continuous surjective map $\Hat\Pi$ on the completed tensor product $\widehat{\otimes}_i$. Since $\Abs0\Omega_c(G)$ is nuclear \cite{as}*{Proposition C.7}, we have $\widehat{\otimes}_i=\widehat{\otimes}_\pi$ \cite{treves}*{Theorem 50.1}, the latter denoting the completed projective tensor product. Since $G$ acts weakly smoothly on $\Abs0\Omega_c(G)$, it acts weakly smoothly on $\Abs0\Omega_c(G)\mathop{\widehat{\otimes}_\pi} E$. Since $E$ is as a Fr\'echet $G_0$- and $\ger g$-representation a quotient of this space, it follows that it is a weakly smooth $G$-representation. 
\end{proof}

\begin{Rem}
  The graded version of the Dixmier--Malliavin theorem offered above (part \eqref{dm-i} of \thmref{Prop}{dm}) admits an independent proof, which does not appeal to \thmref{Cor}{ber-dens}, but rather follows a similar path as Dixmier and Malliavin in their original proof, reducing the statement to low-dimensional cases. To simplify the exposition, we restrict ourselves to the case of a weakly smooth $G$-representation.

  We need to show that $E\subseteq\Pi\Parens1{\Abs0\Omega_c(G)}E$. To that end, we introduce the following terminology: A closed Lie subsupergroup $H$ of $G$ is called \Define{singly generated} if its Lie superalgebra $\ger h$ is generated by a single homogeneous element. 

  Then the following is straightforward: Any singly generated Lie subsupergroup is locally isomorphic to one of the Abelian supergroups $\aff^1$ and $\aff^{0|1}$, or to $\aff^{1|1}$, where the Lie superalgebra has the unique non-zero homogeneous relation $x=[y,y]$. Moreover, there exist singly generated closed Lie subsupergroups $H_1,\dotsc,H_n$ \scth the $n$-fold multiplication morphism $m:H_1\times\dotsm\times H_n\longrightarrow G$ is an isomorphism in a neighbourhood $U$ of the identity. 
  
  Now, fix $v\in E$. We claim the following: For any singly generated sub-supergroup $H$ and any neighbourhood $V\subseteq U$ of $1$, there exist $\omega_0,\omega_1\in\Abs0\Omega_c(H)\subseteq\sh E'(H)\subseteq\sh E'(G)$ and $w\in E$ with $\supp \omega_j\subseteq H\cap V$, \scth $v=\Pi(\omega_0)v+\Pi(\omega_1)w$. Since $\aff^1\subseteq\aff^{1|1}$ as a closed Lie subsupergroup, this follows from Dixmier--Malliavin \cite{dix-mal}*{Theorem 3.3} in case $H_0$ is locally isomorphic to $\aff^1$. 
  
  In case $H$ is isomorphic to $\aff^{0|1}$, we have $\Gamma(\sh O_H)=\knums[\tau]$ where $\tau$ is odd. It follows that $\Abs0{D\tau}$, defined by $\int_H\Abs0{D\tau}\,f=\frac d{d\tau}f$, is a smooth density, and $\int_H\Abs0{D\tau}\,(\tau f)=f(0)$. Thus, the Dirac delta $\delta=\Abs0{D\tau}\,\tau\in\Abs0\Omega_c(H)$ is a smooth density and hence, the statement is obvious in this case. 
  
  Applying the statement inductively, we find $f_0^j,f_1^j\in\Abs0\Omega_c(H_j)\subseteq\sh E'(G)$ and $w_{i_1,\dotsc,i_n}\in E$, $i_j=0,1$, \scth 
  \[
    v=\sum\nolimits_{i_1,\dotsc,i_n=0,1}\Pi(f^1_{i_1}*\dotsm*f_{i_n}^n)w_{i_1,\dotsc,i_n}\ .
  \]
  Now, for $\omega_j\in\Abs0\Omega_c(H_j)$ and $\vphi\in\Gamma(\sh O_G)$, we have 
  \[
    \Dual0{\omega_1*\dotsm*\omega_n}\vphi=\Dual0{\omega_1\otimes\dotsm\otimes\omega_n}{m^\sharp\vphi}.
  \]
  Since $\omega_1\otimes\dotsm\otimes\omega_n$ is in $\Abs0\Omega_c(H_1\times\dotsm\times H_n)$, we find $\omega_1*\dotsm*\omega_n\in\Abs0\Omega_c(G)$, provided that the $\supp\omega_j$ are small enough. This finally proves the claim. 
\end{Rem}

\section{\texorpdfstring{$SF$-representations}{SF-representations}}

In this section, we extend the notion of smooth representations of moderate growth, or \emph{$SF$-representations}, to the case of Lie supergroups. We construct a superalgebra of Schwartz--Berezin densities and show that its representations are in one-to-one correspondence with $SF$-representations of $G$.

\subsection{Schwartz--Berezin densities}

Following \cite{bk}*{2.1}, we will call a measurable function $s:G_0\longrightarrow(0,\infty)$ a \Define{scale} if $s$ and $1/s$ are locally bounded and 
\[
  s(gh)\sle s(g)s(h).
\]
We write $s\preceq s'$ for scales $s,s'$ if is a constant $C>0$ and an integer $N\sge0$ with 
\[
  s(g)\sle Cs'(g)^N
\]
\fa $g\in G_0$. This defines a preorder. The equivalence classes for the largest equivalence relation contained in $\preceq$ are denoted by $[s]$ and called \Define{scale structures}.

\medskip\noindent
In what follows, we fix a scale $s$ on $G_0$. We will always make the assumption that \emph{$s$ dominates the $\ger g$-adjoint scale}, \ie $s\succeq s_\ger g$ where
\[
  s_\ger g(g)\defi\max\Parens1{\Norm0{\Ad(g)|_{\ger g}},\Norm0{\Ad(g^{-1})|_{\ger g}}} 
\]
where we fix some norm on $\ger g$. Observe that there is a constant $C>0$ \scth 
\[
  \max\Parens1{\Abs0{\Delta_\odd(g)},\Abs0{\Delta_\odd(g)}^{-1}}\sle Cs_\ger g(g)^N,
\]
where $N=\dim\ger g_\odd$.

\begin{Def}[sw-dens][Schwartz--Berezin densities]
  We define the space of \Define{Schwartz--Berezin densities} to be 
  \[
    \Sw0{G,[s]}\defi\Set3{\Abs0{Dg}\,f}{\forall u,v\in\Uenv0{\ger g},N\sge0\,:\,\int_{G_0}\Abs0{dg}\,s(g)^N\Abs1{(L_uR_vf)(g)}<\infty},
  \]
  where $\Abs0{Dg}$ and $\Abs0{dg}$ are some choices of left invariant Berezin density on $G$ resp.~left invariant density on $G_0$. This space is endowed with the locally convex topology generated by the seminorms
  \[
    p^s_{u,v,N}(\Abs0{Dg}\,f)\defi\int_{G_0}\Abs0{dg}\,s(g)^N\Abs1{(L_uR_vf)(g)}.
  \]
  Clearly, the locally convex super-vector space $\Sw0{G,[s]}$ is independent of the choice of $\Abs0{Dg}$, $\Abs0{dg}$, and the representative $s$ of the scale structure $[s]$.
\end{Def}

Similarly, there is a space of Schwartz densities $\Sw0{G_0,[s]}\subseteq\Abs0\Omega(G_0)$. According to \cite{bk}*{2.5}, it is defined as the space of smooth vectors for the bi-regular representation $L_0\times R_0$ of $G_0\times G_0$ on the space $\sh R(G_0,[s])$, the set of \emph{continuous} densities $\omega$ that are \Define{rapidly decreasing} in the sense that
\[
  \forall N\in\nats\,:\,\int_{G_0}\Abs0{\omega} s^N<\infty.
\]
We have the following description of $\Sw0{G,[s]}$ in terms of $\Sw0{G_0,[s]}$.

\begin{Prop}[sw-iso]
  The isomorphism from \thmref{Cor}{ber-dens} induces an isomorphism 
  \begin{equation}\label{eq:sw-iso}
    \Sw0{G,[s]}=\Uenv0{\ger g}\otimes_{\Uenv0{\ger g_\ev}}\Sw0{G_0,[s]}.    
  \end{equation}
  In particular, $\Sw0{G,[s]}$ is nuclear space and $G$-invariant for the left regular representation $L$, as a subspace of $\Abs0\Omega(G)$.
\end{Prop}

\begin{proof}
  By the above definitions, we have $\Abs0{Dg}\,f\in\Sw0{G,[s]}$ if and only if for any $u,v\in\Uenv0{\ger g}$, we have 
  \[
    \omega\defi\Abs0{dg}\,j^\sharp(L_uR_vf)\in\sh R(G_0,[s]),
  \]
  where we abbreviate $j\defi j_{G_0}$. Such a density $\omega$ is smooth, and for $x\in\ger g_\ev$, we have 
  \[
    L_x\omega=\Abs0{dg}\,j^\sharp(L_{xu}R_vf)\in\sh R(G_0,[s]).
  \]
  One argues similarly for $R_x$, so that $\omega\in\Sw0{G_0,[s]}$.

  Now, let $\Abs0{Dg}\,f\in\Abs0\Omega(G)$. We may assume w.l.o.g.~that $\Abs0{Dg}$ and $\Abs0{dg}$ are related by Equation \eqref{eq:invber}. Then Equation \eqref{eq:invber-haar} implies that $\Abs0{Dg}\,f$ corresponds to 
  \[
    \sum\nolimits_i\gamma_i'\otimes\Abs0{dg}\,\Delta_\odd\,j^\sharp\Parens1{L_{S(\gamma_i'')}(f)}
  \]
  where $\Delta(\gamma)=\sum_i\gamma_i'\otimes\gamma_i''$. For any $u\in\Uenv0{\ger g}$, we therefore have 
  \[
    \Abs0{dg}\,\Delta_\odd\,j^\sharp\Parens1{L_u(f)}\in\Sw0{G_0,[s]},
  \]
  since $s$ dominates the $\ger g$-adjoint scale by assumption.

  Conversely, let the Berezinian density $\Abs0{Dg}\,f$ correspond to $u\otimes\omega$, where we assume $\omega=\Abs0{dg}\,h\in\Sw0{G_0,[s]}$ and $u\in\Uenv0{\ger g}$. By \thmref{Cor}{ber-dens}, we have $f=L_u(\psi h)$, with $\psi$ defined in Equation \eqref{eq:psidef}. 

  For $v,w\in\Uenv0{\ger g}$, $g\in G_0$, we expand 
  \[
    (-1)^{\Abs0u\Abs0w}\Delta(vu)=\sum\nolimits_iv_i'\otimes v_i'',\quad
    \Delta(w)=\sum\nolimits_jw_j'\otimes w_j''.
  \]
  Then we compute for $\vkappa_{ij}\defi(-1)^{\Abs0\psi(\Abs0{v_i''}+\Abs0{w_j''})+\Abs0{v_i'}\Abs0{w_j'}}$ that
  \begin{equation}\label{eq:psicomp}
    j^\sharp(L_vR_w(\psi f))(g)=\sum\nolimits_{i,j}\vkappa_{ij}\Parens1{L_{v_i'}R_{w_j'}(\psi)}(g)\Parens1{L_{v_i''}R_{w_j''}(h)}(g)
  \end{equation}
  with 
  \begin{align*}
    \Parens1{L_{v'}R_{w'}(\psi)}(g)&=\Parens1{R_{\iota(\Ad(g^{-1})(S(v'))w')}\Delta_\odd^{-1}}(g)\\
    &=\delta_\odd\Parens1{S(\iota(\Ad(g^{-1})(S(v'))w'))}\Delta_\odd(g)^{-1}.    
  \end{align*}
  We have 
  \[
    \iota(xay)=x\iota(a)(y-\delta_\odd(y))
  \]
  \fa $x,y\in\ger g_\ev$ and $a\in\Uenv0{\ger g}$ \cite{gorelik-ghost}*{Equation (3)}. Moreover, $\delta_\odd$ is a character of $\Uenv0{\ger g_\ev}$ and in particular $\Ad(G_0)$-invariant. Finally, there is a constant $C>0$ \scth 
  \[
    \Abs1{\iota\Parens1{\beta(\Ad(g)(\xi)\eta}}=\Abs3{\int_{\ger g_\odd}\Ad(g)(\xi)\eta}\sle C\Norm0{\Ad(g)|_{\ger g_\odd}}^k\Norm0\xi\Norm0\eta
  \]
  \fa $\xi\in\bigwedge^k\ger g_\odd$, $\eta\in\bigwedge\ger g_\odd$. (Here, $\Norm0\cdot$ denotes some submultiplicative norm on $\bigwedge\ger g_\odd$.) It follows that there exist a constant $C>0$ and an integer $N\sge 0$ \scth
  \[
    \Abs1{\delta_\odd\Parens1{S(\iota(\Ad(g^{-1})(S(v'))w'))}}\sle Cs_\ger g(g)^N.
  \]
  for all $g\in G_0$. The sum in Equation \eqref{eq:psicomp} is finite, so we may conclude that $\Abs0{dg}\,j^\sharp(L_vR_w(f))$ is a Schwartz density on $(G_0,[s])$ once so is $\Abs0{dg}\,\Parens1{L_{u''}R_{v''}(h)}$.

  To that end, similarly as above, we note that 
  \[
    j^\sharp\Parens1{L_{xa}R_{yb}(h)}=L_xR_yj^\sharp\Parens1{L_aR_bh}
  \]
  \fa $x,y\in\ger g_\ev$, $a,b\in\Uenv0{\ger g}$, that 
  \[
      \Parens1{R_{\beta(\Ad(g)(\xi)\eta)}h}(g)=\eps(\beta(\Ad(g)(\xi)\eta)),
  \]
  and that there is a constant $C>0$ \scth
  \[
    \Abs1{\eps(\beta(\Ad(g)(\xi)\eta))}\sle C\Norm0{\Ad(g)|_{\ger g_\odd}}^k\Norm0\xi\Norm0\eta.
  \]
  Thus, there is an integer $N\sge 0$ \scth for all $g\in G_0$, we have 
  \[
    \Abs1{\Parens1{L_{v''}R_{w''}(h)}(g)}\sle s_\ger g(g)^N\Abs0{H(g)}
  \]
  where $H=\sum_\ell(L_{a_\ell}R_{b_\ell}h)$ \fs $a_\ell,b_\ell\in\Uenv0{\ger g_\ev}$. 

  In summary, we have shown the isomorphism in Equation \eqref{eq:sw-iso}, in particular, $\Sw0{G,[s]}$ is a $G$-invariant subspace of $\Abs0\Omega(G)$. Inspecting the above formul\ae{}, it is evident that it is an isomorphism of topological vector spaces, if $\Uenv0{\ger g}\otimes_{\Uenv0{\ger g_\ev}}\Sw0{G_0,[s]}$ is endowed with the natural topology on $\bigwedge\ger g_\odd\otimes\Sw0{G_0,[s]}$. The nuclearity now follows from that of $\Sw0{G_0,[s]}$ \cite{bk}*{Corollary 5.6}.
\end{proof}

\begin{Prop}[sw-conv]
  The subspace $\Sw0{G,[s]}\subseteq\sh D'(G)$ is bi-invariant under the regular representation of $G$. Via the isomorphism in Equation \eqref{eq:sw-iso}, it inherits a non-unital Fr\'echet superalgebra structure with continuous multiplication, determined uniquely by the following facts:
  \begin{enumerate}[wide]
    \item\label{sw-conv-i} The following is a non-unital graded subalgebra bi-invariant under $G_0$:
    \[
      \Sw0{G_0,[s]}=\Uenv0{\ger g_\ev}\otimes_{\Uenv0{\ger g_\ev}}\Sw0{G_0,[s]}.
    \]
    \item\label{sw-conv-ii} For $u,v\in\Uenv0{\ger g}$ and $\omega\in\Sw0{G_0,[s]}$, we have 
    \[
      u*(v\otimes\omega)=(u\otimes1)*(v\otimes\omega)=uv\otimes\omega.
    \]
    \item\label{sw-conv-iii} For $\omega\in\Sw0{G_0,[s]}$ and $u\in\Uenv0{\ger g}$, the product $u*\omega$ is given by 
    \begin{equation}\label{eq:sw-conv}
        \int_{G_0}(u*\omega)\vphi=\int_{G_0}\omega\,j_{G_0}^\sharp\Parens1{L_{S(u)}\vphi}=\Dual1{u\otimes\omega}\vphi,
    \end{equation}
    \fa compactly supported superfunctions $\vphi\in\Gamma_c(\sh O_G)$.
  \end{enumerate}
\end{Prop}

\begin{proof}
  We already know that $\Sw0{G,[s]}$ is invariant under $L_u$ and $L_g$ for any $u\in\Uenv0{\ger g}$, $g\in G_0$. To see that $\Sw0{G,[s]}$ is also invariant under the right regular action $R$, it will be sufficient to show that $\Sw0{G,[s]}$ is stable under $(-)^\vee$, defined in Equation \eqref{eq:checkdef}.

  Choose a basis $x_1,\dotsc,x_q$ of $\ger g_\odd$, and let $x^1,\dotsc,x^q$ be the dual basis of $\ger g_\odd^*$. We write $x_I\defi x_{i_1}\dotsm x_{i_k}\in\bigwedge\ger g_\odd$ and $x^I\defi x^{i_1}\dotsm x^{i_k}\in\bigwedge\ger g_\odd^*$ for $I=(1\sle i_1<\dotsm<i_k\sle q)$. Then we compute by Equation \eqref{eq:tensor-dist}
  \begin{align*}
    \Dual1{(\beta(x_I)\otimes\omega)^\vee}\vphi&=(-1)^{\Abs0I\Abs0\vphi}\int_{G_0}\omega(g)\vphi(S(\beta(x_I));g^{-1})\\
    &=\sum\nolimits_J\int_{G_0}(-1)^{\Abs0J\Abs0\vphi}\check\omega(g)\Dual0{x^J}{\Ad(g)(x_I)}\vphi(\Ad(g^{-1})(S(x_J);g)\\
    &=\Dual2{\sum\nolimits_J\beta(x_J)\otimes\check\omega\Dual0{x^J}{\Ad(\cdot)(x_I)}}\vphi
  \end{align*}
  \fa $\vphi\in\Gamma_c(\sh O_G)$. Here, observe that $\Abs0J=\Abs0I$, because the adjoint action by $G_0$ on $\bigwedge\ger g_\odd$ respects the $\ints$-grading. 

  By the assumption on the scale $s$, we have $\check\omega\in\Sw0{G_0,[s]}$ and thus 
  \[
    \Abs0{dg}\,\sum\nolimits_J\check\omega\Dual0{x^J}{\Ad(\cdot)(x_I)}\in\Sw0{G_0,[s]}.
  \]
  In view of \thmref{Prop}{sw-iso}, this shows that $(\beta(x_I)\otimes\omega)^\vee\in\Sw0{G,[s]}$. Therefore, $\Sw0{G,[s]}$ is invariant under $(-)^\vee$ and bi-invariant under $G$.

  In follows that there is a well-defined operation $*$ on $\Sw0{G,[s]}$, defined by 
  \begin{equation}\label{eq:sw-convdef}
    (u\otimes\omega)*(v\otimes\varpi)\defi\sum\nolimits_juv_j\otimes(\omega_j*\varpi),
  \end{equation}
  for arbitrary $u,v\in\Uenv0{\ger g}$ and $\omega,\varpi$, where we decompose $R_{S(v)}\omega=\sum_jv_j\otimes\omega_j$. 

  If $\omega,\varpi$ are compactly supported, then by Equation \eqref{eq:tensor-dist}, we have 
  \begin{align*}
    \bigl\langle(u\otimes\omega)&*(v\otimes\varpi),\vphi\bigr\rangle
    =\sum\nolimits_j(-1)^{\Abs0{uv_j}\Abs0\vphi}\int_{G_0}(\omega_j*\varpi)(g)\vphi\Parens1{\Ad(g^{-1})(uv_j);g}\\
    &=\sum\nolimits_j(-1)^{\Abs0{uv_j}\Abs0\vphi}\int_{G_0\times G_0}\omega_j(g)\varpi(h)\vphi\Parens1{\Ad((gh)^{-1})(uv_j);gh}\\
    &=\sum\nolimits_j(-1)^{\Abs0{uv_j}\Abs0\vphi}\int_{G_0\times G_0}\omega_j(g)\varpi(h)m^\sharp(\vphi)\Parens1{\Ad(g^{-1})(uv_j)\otimes1;g,h}\\
    &=\Dual2{\sum\nolimits_j(uv_j\otimes\omega_j)\otimes\varpi}{m^\sharp(\vphi)}\\
    &=\Dual1{L_uR_{S(v)}\omega\otimes\varpi}{m^\sharp(\vphi)}\\
    &=(-1)^{\Abs0{uv}\Abs0\vphi}\int_{G_0\times G_0}\omega(g)\varpi(h)m^\sharp(\vphi)(\Ad(g^{-1})(u)v\otimes1;g,h)\\
    &=(-1)^{\Abs0{uv}\Abs0\vphi}\int_{G_0\times G_0}\omega(g)\varpi(h)m^\sharp(\vphi)(\Ad(g^{-1})(u)\otimes\Ad(h^{-1})(v);g,h)\\
    &=\Dual1{(u\otimes\omega)\otimes(v\otimes\varpi)}{m^\sharp(\vphi)},
  \end{align*}
  so that $*$ on $\Sw0{G,[s]}$ extends the convolution on $\Abs0\Omega_c(G)$. By \thmref{Prop}{sw-iso}, $\Abs0\Omega_c(G)\subseteq\Sw0{G,[s]}$ is a dense subspace and $\Sw0{G,[s]}$ is nuclear. To finish the proof of our assertion, it is by \cite{treves}*{Theorem 50.1} sufficient to show that the convolution $*$ on $\Abs0\Omega_c(G)$ is separately continuous in the topology induced by $\Sw0{G,[s]}$.

  Since $(-)^\vee$ is continuous on $\Sw0{G,[s]}$, it will be sufficient to show continuity in the second argument. In view of \thmref{Cor}{dens-conv} \eqref{dens-conv-i}--\eqref{dens-conv-iii}, we have the identity 
  \[
    L_uR_v(\omega*\varpi)=(-1)^{\Abs0v\Abs0\omega}(L_u\omega)*(R_v\varpi)
  \]
  \fa $u,v\in\Uenv0{\ger g}$ and $\omega,\varpi\in\Abs0\Omega_c(G)$. Together with the fact that for any $v\in\Uenv0{\ger g}$, $R_v$ is continuous on $\Sw0{G,[s]}$, it follows that it is sufficient to show that 
  \[
    \varpi\longmapsto p_{1,1,N}^s\Parens1{(L_uR_{S(v)}\omega)*\varpi}:\Sw0{G_0,[s]}\longrightarrow\reals
  \]
  is a continuous seminorm for any integer $N\sge0$. But this follows from the continuity of the convolution on $\Sw0{G_0,[s]}$ \cite{wallach1}*{Theorem 7.1}.
\end{proof}

\subsection{\texorpdfstring{$SF$-representations of supergroups}{SF-representations of supergroups}}

Fix a scale $s$ on $G_0$ dominating the $\ger g$-adjoint scale. Recall \cite{bk}*{Definition 2.6, Lemma 2.10} that a continuous Fr\'echet $G$-representation $\pi_0$ on $E$ is called an \Define{$F$-representation} or a Fr\'echet representation of \Define{moderate growth} of $(G_0,[s])$ if the topology of $E$ is generated by a countable collection $(p_j)$ of seminorms \scth for any $j$, there exist an index $k$, a constant $C>0$, and an integer $N\sge0$ with 
\begin{equation}\label{eq:f}
  p_j\Parens1{\pi_0(g)v}\sle Cs(g)^Np_k(v)
\end{equation}
\fa $v\in E$ and $g\in G_0$. It is called an \Define{$SF$-representation} or \Define{smooth} if it is in addition weakly smooth.

\medskip\noindent
In view of this terminology, we make the following definition.

\begin{Def}[sf-rep][\protect{$SF$-representations}]
  Let $\pi$ be a continuous representation $G$ on a Fr\'echet super-vector space $E$. Then $\pi$ is called an \Define{$F$-representation} of $(G,[s])$ if the topology of $E$ is generated by a countable collection $(p_j)$ of seminorms \scth for some norm $\Norm0\cdot$ on $\bigwedge\ger g_\odd$ and for any index $j$, there is an index $k$, a constant $C>0$, and an integer $N\sge0$ with 
  \begin{equation}\label{eq:super-f}
    p_j\Parens1{d\pi(\beta(\eta))\pi_0(g)v}\sle C\Norm0{\eta}s(g)^Np_k(v)
  \end{equation}
  \fa $v\in E_\infty$, $g\in G_0$, and $\eta\in\bigwedge\ger g_\odd$. If in addition, $\pi$ is a weakly smooth $G$-representation, then it is called an \Define{$SF$-representation} of $(G,[s])$.
\end{Def}

\begin{Rem}
  If $\pi$ is an $F$-representation (resp.~an $SF$-representation) of $(G,[s])$, then $\pi_0$ is an $F$-representation (resp.~an $SF$-representation) of $(G_0,[s])$. Indeed, $E_\infty$ is dense in $E$, and taking $\eta=1$ in Equation \eqref{eq:super-f}, we obtain Equation \eqref{eq:f}. Also by definition, if $\pi$ is an $F$-representation of $(G,[s])$ on $E$, then the subrepresentation on the space $E_\infty$ of smooth vectors is an $SF$-representation \cite{bk}*{Corollary 2.16}. 

  In particular, using \cite{bk}*{(2.2)}, \thmref{Prop}{sw-iso}, and \thmref{Prop}{sw-conv}, we find that the left and right regular representations $L$ and $R$ on $\Sw0{G,[s]}$ are $SF$-representations of $(G,[s])$.  
\end{Rem}

In fact, the $F$-representations are characterised among the continuous representations of $G$ by the growth of the underlying $G_0$-representation.

\begin{Lem}[f-char]
  Let $\pi$ be a continuous (resp.~weakly smooth) representation of $G$ on a Fr\'echet super vector-space $E$. Then $\pi$ is an $F$-representation (resp.~an $SF$-representation) of $(G,[s])$ if and only if $\pi_0$ is an $F$-representation (resp.~an $SF$-representation) of $(G_0,[s])$.
\end{Lem}

\begin{proof}
  It is sufficient to consider the case of $F$-representations. As noted above, if $\pi$ is an $F$-representation of $(G,[s])$, then $\pi_0$ is an $F$-representation of $(G_0,[s])$. Conversely, assume that $\pi_0$ is an $F$-representation of $(G_0,[s])$. Since $s$ dominates the $\ger g$-adjoint scale, the adjoint representation of $G_0$ on $\bigwedge\ger g_\odd$ is an $F$-representation. Hence, so is $\bigwedge\ger g_\odd\otimes E$. Manifestly, this gives the condition in Equation \eqref{eq:super-f}.
\end{proof}

\begin{Rem}
  From \thmref{Lem}{f-char}, we obtain the following: Let $\pi$ be a continuous $G$-representation on a Banach super vector-space $E$. Then $\pi$ is an $F$-representation of $(G,[s])$ if and only if $\pi_0$ is $s$-bounded in the sense that $s\succeq s_{\pi_0}$ where 
  \[
    s_{\pi_0}(g)\defi\max\Parens1{\Norm0{\pi_0(g)},\Norm0{\pi_0(g^{-1})}}.
  \]
  In particular, in this case, the $G$-representation on $E_\infty$ is an $SF$-representation \cite{bk}*{Corollary 2.16}.
\end{Rem}

\medskip\noindent
For $F$-representations of $G$, we obtain the following variant of the Dixmier--Malliavin theorem, generalising \cite{bk}*{Proposition 2.20}. Compare \cite{ducloux}*{Exemple 2.3.3}.

\begin{Prop}[sw-dm]
  Let $E$ be a Fr\'echet super-vector space over $\knums$. Then we have the following facts: 
  \begin{enumerate}[wide]
    \item\label{swdm-i} If $E$ carries the structure of an $F$-representation $\pi$ of $(G,[s])$, then the integrated action $\Pi$ extends continuously to an action of $\Sw0{G,[s]}$, also called the \Define{integrated action} of $\pi$. We have the equality
    \begin{equation}\label{eq:sw-dm}
      E_\infty=\Pi\Parens1{\Sw0{G,[s]}}E=\Pi\Parens1{\Sw0{G,[s]}}E_\infty.
    \end{equation}
    \item\label{swdm-ii} Conversely, let $\Sw0{G,[s]}$ act continuously and non-degenerately \via $\Pi$ on $E$. Then $\Pi$ is integrated from a unique $SF$-representation of $(G,[s])$.
  \end{enumerate}
  In particular, we obtain an equivalence of the category of $SF$-representations of $(G,[s])$ with the category of non-degenerate continuous Fr\'echet $\Sw0{G,[s]}$-modules. 
\end{Prop}

\begin{proof}
  If $E$ is an $F$-representation of $(G_0,[s])$, then $\Sw0{G_0,[s]}$ acts continuously on $E$, and $E_\infty=\Sw0{G_0,[s]}E=\Sw0{G_0,[s]}E_\infty$ \cite{bk}*{Proposition 2.20}. Conversely, if $E$ carries a continuous non-degenerate action of $\Sw0{G_0,[s]}$, then this action is integrated from a unique $SF$-representation of $(G_0,[s])$ (\loccit). Using these facts, together with \thmref{Prop}{sw-iso} and \thmref{Prop}{sw-conv}, the proof of the claim is the same as that of \thmref{Prop}{dm}. We therefore leave the details to the reader. 
\end{proof}

\section{Harish-Chandra supermodules}\label{sec:cw}

In this section, we come to our main result, a generalisation of the Casselman--Wallach theorem to supergroups. 

\subsection{Basic facts and definitions}

In what follows, we assume that the underlying Lie group $G_0$ of $G$ is almost connected and real reductive \cite{wallach1}*{2.1} and let $K_0\subseteq G_0$ be a maximal compact subgroup. We fix on $G_0$ the maximal scale structure \cite{bk}*{2.1.1} and omit the mention of $[s]$ in our notation. In particular, any Banach representation of $G$ is an $F$-representation.

\begin{Def}[hc-def][Harish-Chandra supermodules]
  A \Define{$(\ger g,K_0)$-module} is by definition a complex, $\ints/2\ints$ graded, locally finite $K_0$-representation $V$, endowed with a $K_0$-equivariant $\ger g$-module structure, which extends the derived $\ger k_0$-action on $V$. A \Define{morphism of $(\ger g,K_0)$-modules} $\phi:U\longrightarrow V$ is an even $\cplxs$-linear map that is equivariant for the actions of $\ger g$ and $K_0$.

  A $(\ger g,K_0)$-module is called \Define{Harish-Chandra} or a \Define{Harish-Chandra supermodule} if it is $K_0$-multiplicity finite and finitely generated over $\Uenv0{\ger g}$. The full subcategory of the category of $(\ger g,K_0)$-modules whose objects are the Harish-Chandra supermodules is denoted by $\HC(\ger g,K_0)$.
\end{Def}

The following observation is elementary, but effective.

\begin{Lem}[hc-char]
  Let $V$ be a $(\ger g,K_0)$-module. Then $V\in\HC(\ger g,K_0)$ if and only if its restriction $V|_{(\ger g_\ev,K_0)}$ to a $(\ger g_\ev,K_0)$-module lies in $\HC(\ger g_\ev,K_0)$.
\end{Lem}

\begin{proof}
  We need only observe that $\Uenv0{\ger g}$ is finitely generated as a $\Uenv0{\ger g_\ev}$-module.
\end{proof}

\begin{Lem}[kfin]
  Let $E$ be an $F$-representation of $G$ (for instance, a Banach representation). Then the space $\smash{E_\infty^{(K_0)}}$ of $K_0$-finite and smooth vectors is a $(\ger g,K_0)$-module.
\end{Lem}

\begin{proof}
  Since the action of $K_0$ on $\Uenv0{\ger g}$ is locally finite, we see that the $\ger g$-action on $E_\infty$ leaves $\smash{E^{(K_0)}_\infty}$ invariant. 
\end{proof}

\begin{Rem}
  Let $\pi_0$ be a continuous $G$-representation on a complex Banach super-vector space $E$. Denoting by $C$ the Casimir element of $\ger g_\ev$, assume that either 
  \begin{enumerate}
    \item $d\pi_0(C)\in\End0{E_\infty}$ extends continuously to $E$, or 
    \item $P(d\pi_0(C))=0$ on $E_\infty$ for some polynomial $P$.
  \end{enumerate}
  Then it is known that the space $\smash{E^{(K_0)}}$ of $K_0$-finite vectors is contained in $E_\infty$ \cite{bk}*{Corollary 3.10}. Hence, if $\pi_0$ is the $G_0$ part of a continuous $G$-representation, then $\smash{E^{(K_0)}}$ is a $(\ger g,K_0)$-module, by \thmref{Lem}{kfin}.
\end{Rem}

\subsection{Globalisation of Harish-Chandra supermodules} 

\begin{Def}[cw][Casselman--Wallach representations]
  An $SF$-representation $(E,\pi)$ of $G$ is called \Define{Casselman--Wallach} or a \Define{CW representation} if the space $\smash{E^{(K_0)}}$ of $K_0$-finite and smooth vectors is in $\HC(\ger g,K_0)$.

  If $V\in\HC(\ger g,K_0)$, then an isomorphism $\phi:V\longrightarrow \smash{E^{(K_0)}}$ of $(\ger g,K_0)$-modules, where $(E,\pi)$ is an $SF$-representation, is called an \Define{$SF$-globalisation} of $V$. Any $SF$-globalisation of a Harish-Chandra supermodule is a CW representation of $G$.

  A CW globalisation $\phi:V\longrightarrow E$ is called \Define{minimal} if for any CW globalisation $\psi:V\longrightarrow H$, there exists an even continuous $G$-equivariant map $\tilde\psi:E\longrightarrow H$ \scth $\tilde\psi\circ\vphi=\vphi$. Since the $K_0$-finite vectors are dense in $E$, such a map $\tilde\psi$ is unique. Thus, minimal globalisations (if they exist) are unique up to canonical isomorphism.

  Dually, a CW globalisation $\phi:V\longrightarrow E$ is called \Define{maximal} if for any CW globalisation $\psi:V\longrightarrow H$, there exists an even continuous $G$-equivariant map $\tilde\psi:H\longrightarrow E$ \scth $\tilde\psi\circ\psi=\vphi$. Again, maximal globalisations (if they exist) are unique up to canonical isomorphism.
\end{Def}

We are now ready to state our main theorem. 

\begin{Th}[cw][super Casselman--Wallach theorem]
  Let $V\in\HC(\ger g,K_0)$. Up to isomorphism, there is a unique CW globalisation of $V$. 
\end{Th}

We postpone the proof to Subsection \ref{sub:proof} and give a number of corollaries. The derivation of these follows the same procedures as in the Lie group case \cite{wallach2}*{11.6.8}.

\begin{Cor}
  The functor mapping $(E,\pi)$ to $E_\infty^{(K_0)}$ sets up an additive equivalence between the category $\CW(G)$ of CW representations of $G$ and the category $\HC(\ger g,K_0)$ of Harish-Chandra supermodules. In particular, the category $\CW(G)$ is Abelian.
\end{Cor}

\begin{Cor}[topmor]
  Let $f:E\longrightarrow F$ be a morphism of CW $G$-representations. Then $f$ is a topological morphism with closed image. 

  Here, $f:E\longrightarrow F$ is called a \Define{topological morphism} if the induced map
  \[
     E/\ker f\longrightarrow\im f
  \] 
  is an isomorphism of topological vector spaces. 
\end{Cor}

As a corollary to the proof of \thmref{Th}{cw}, we obtain the following.

\begin{Cor}
  Any $E\in\CW(G)$ is the space of smooth vectors of a continuous Hilbert $G$-representation. 
\end{Cor}

\subsection{Proof of Theorem \ref{Th:cw}}\label{sub:proof}

Having stated our main result, together with some immediate corollaries, let us come to its proof. 

\begin{proof}[\prfof{Th}{cw}]
  First, we show that $V$ has a minimal $SF$-globalisation $V_+\supseteq V$. We mimic the construction detailed in \cite{bk}*{\S~6}. 

  By \thmref{Lem}{hc-char}, we have $U\defi V|_{(\ger g_\ev,K_0)}\in\HC(\ger g_\ev,K_0)$. Thus, there is a finite set $v_1,\dotsc,v_n$ of homogeneous vectors generating the $\Uenv0{\ger g_\ev}$-module $V$ and a continuous Hilbert representation $(E,\pi_0)$ of $G_0$ \scth $\smash{E^{(K_0)}_\infty}=U$ \cite{bk}*{\S~5.1}.

  Since $\Sw0G$ is invariant under $(-)^\vee$, \thmref{Prop}{sw-iso} shows that the map
  \begin{equation}\label{eq:sw-r-iso}
    \Sw0{G_0}\otimes_{\Uenv0{\ger g_\ev}}\Uenv0{\ger g}\longrightarrow\Sw0G:\omega\otimes u\longmapsto R_{S(u)}(\omega)
  \end{equation}
  is an isomorphism of right $\Uenv0{\ger g}$-modules. Here, $\Uenv0{\ger g}$ acts from the right on $\Sw0G$ by $\omega u\defi(-1)^{\Abs0\omega\Abs0u}R_{S(u)}(\omega)$. We define, for $\omega\in\Sw0{G_0}$ and $v\in V$
  \begin{equation}\label{eq:swact-def}    
    \Pi(\omega)v\defi\sum\nolimits_i\Pi_0(\omega_i)u_iv,
  \end{equation}
  where 
  \[
    \omega=\sum\nolimits_IR_{S(u_i)}(\omega_i)
  \]
  is any decomposition with $\omega_i\in\Sw0{G_0}$ and $u_i\in\Uenv0{\ger g}$. To see that this is well-defined, we need only remark that 
  \[
    \Pi_0(R_{-x}\omega)=\Pi_0(\omega)d\pi_0(x)
  \]
  \fa $\omega\in\Sw0{G_0}$ and $x\in\ger g_\ev$.

  Now, consider the graded subspace $\mathscr N\subseteq\Sw0G^n$, defined by 
  \[
    \mathscr N\defi\Set2{(\omega_1,\dotsc,\omega_n)\in\Sw0G^n}{\sum\nolimits_j\Pi(\omega_j)v_j=0}.
  \]
  We claim that it is closed and invariant under the action of $\Sw0G$ by left convolution. To prepare the proof of this claim, we briefly suspend our argument and establish some ancillary lemmas.
\end{proof}

Let $\tilde V\in\HC(\ger g,K_0)$ be the dual Harish-Chandra module of $V$, defined as the set of $K_0$-finite vectors in the algebraic dual $V^*$. Then $\tilde V$ is also the dual of $\smash{V|_{(\ger g_\ev,K_0)}}$ \cite{bk}*{\S~4}, and in particular $\tilde V$ is contained in the space $\tilde E$ of continuous vectors of the topological dual $E'$ of $E$ \cite{bk}*{Lemma 5.3}.

\begin{Lem}[hcmod-pairid]
  Let $v\in V$, $u\in\Uenv0{\ger g}$ and $g\in G_0$. We have the identity
  \[
    \Dual0\xi{\pi_0(g^{-1})\Ad(g)(u)v}=(-1)^{\Abs0\xi\Abs0u}\Dual0{S(u)\xi}{\pi_0(g^{-1})v}.
  \]
\end{Lem}

\begin{proof}
  The equality is obvious for $g\in K_0$. Since $G_0'K_0=G_0$, where $G_0'$ is the connected component of the identity of $G_0$, we may assume that $G_0$ is connected. 

  To prove the assertion in that case, assume first that $u\in\beta(\bigwedge\ger g_\odd)$. The image $F$ of $\beta(\bigwedge\ger g_\odd)$ in $\End0V$ is finite-dimensional, so the linear map 
  \[
    F\longrightarrow V\subseteq E:u\longmapsto uv
  \]
  is continuous. For $x\in\ger g_0$, we may hence exchange limits and compute 
  \[
    \frac d{dt}\Big|_{t=0}\Ad(\exp(tx))(u)v=[x,u]v=d\pi_0(x)uv-uxv.
  \]
  Thus, we have 
  \[
    \frac d{dt}\Big|_{t=0}\pi_0(\exp(-tx))\Ad(\exp(tx))(u)v=-d\pi_0(x)uv+[x,u]v=-uxv,
  \]
  by the smoothness of the $G_0$-representation $E_\infty$. Hence 
  \begin{align*}
    \frac d{dt}\Big|_{t=0}\Dual0\xi{\pi_0(e^{-tx})\Ad(e^{tx})(u)v}&=-(-1)^{\Abs0\xi\Abs0u}\Dual0{S(u)\xi}{d\pi_0(x)v}\\
    &=(-1)^{\Abs0\xi\Abs0u}\frac d{dt}\Big|_{t=0}\Dual0{S(u)\xi}{\pi_0(e^{-tx})v}.    
  \end{align*}
  By the uniqueness of initial value problems, the equality follows for $g=e^x$. Since $\exp$ is a local diffeomorphism and $G_0$, being connected, is generated by a neighbourhood of the identity, the equality holds for arbitrary $g\in G_0$.

  To remove the restriction on $u$, recall that $\Uenv0{\ger g}=\Uenv0{\ger g_\ev}\beta(\bigwedge\ger g_\odd)$. By linearity in $u$, it is sufficient to consider $u=u'u''$ for $u'\in\Uenv0{\ger g_\ev}$ and $u''\in\beta(\bigwedge\ger g_\odd)$. Then
  \begin{align*}
    \Dual1\xi{\pi_0(g^{-1})\Ad(g)(u)v}&=\Dual1\xi{\pi_0(g^{-1})\Ad(g)(u')\Ad(g)(u'')v}\\
    &=\Dual1{S(u')\xi}{\pi_0(g^{-1})\Ad(g)(u'')v}\\
    &=(-1)^{\Abs0\xi\Abs0u}\Dual1{S(u'')S(u')\xi}{\pi_0(g^{-1})v}\\
    &=(-1)^{\Abs0\xi\Abs0u}\Dual1{S(u)\xi}{\pi_0(g^{-1})v}.    
  \end{align*}
  This proves the claim in general.
\end{proof}

\begin{Lem}[matrixcoeff-l-id]
  For $u\in\Uenv0{\ger g}$, $\omega\in\Sw0G$, $v\in V$, and $\xi\in\tilde V$, we have
  \[
    \Dual1\xi{\Pi(L_u(\omega))v}
    =(-1)^{\Abs0\xi\Abs0u}\Dual1{S(u)\xi}{\Pi(\omega)v}.
  \]
\end{Lem}

\begin{proof}
  For $v\in V$ and $\xi\in\tilde V$, we define $M_{\xi,v}\in\Gamma(\sh O_G)$ by 
  \[
    M_{\xi,v}(u;g)\defi (-1)^{\Abs0u\Abs0v}\Dual0\xi{\pi_0(g)v}.
  \]
  Clearly, this is well-defined. 

  For $u\in\Uenv0{\ger g}$ and $\omega\in\Sw0{G_0}$, we compute 
  \begin{align*}
    \Dual1\xi{\Pi(R_{S(u)}(\omega))v}&=(-1)^{\Abs0u\Abs0v}\Dual1\xi{\Pi_0(\omega)uv}=(-1)^{\Abs0u\Abs0v}\int_{G_0}\omega(g)\Dual0\xi{\pi_0(g)uv}\\
    &=(-1)^{\Abs0\xi\Abs0u}\Dual1{R_{S(u)}(\omega)}{M_{\xi,v}}.    
  \end{align*}
  By Equation \eqref{eq:sw-r-iso}, it follows that 
  \[
    \Dual0\xi{\Pi(\omega)v}=(-1)^{\Abs0\xi\Abs0\omega}\Dual0\omega{M_{\xi,v}}
  \]
  for any $\omega\in\Sw0G$. In particular, if $u\in\Uenv0{\ger g}$, we have 
  \begin{align*}
    \Dual1\xi{\Pi(L_u(\omega))v}
    &=(-1)^{\Abs0\xi(\Abs0u+\Abs0\omega)+\Abs0u\Abs0\omega}\Dual1\omega{L_{S(u)}(M_{\xi,v})}\\
    &=(-1)^{\Abs0\xi(\Abs0u+\Abs0\omega)+\Abs0u\Abs0\omega}\Dual1\omega{M_{S(u)\xi,v}}
    =(-1)^{\Abs0\xi\Abs0u}\Dual1{S(u)\xi}{\Pi(\omega)v},    
  \end{align*}
  since 
  \begin{align*}
    L_u(M_{\xi,v})(u';g)&=(-1)^{\Abs0u\Abs0\xi+\Abs0{u'}\Abs0v}\Dual1\xi{\pi_0(g)\Ad(g^{-1})(S(u))u'v}\\
    &=(-1)^{\Abs0{u'}\Abs0v}\Dual1{uv}{\pi_0(g)u'v}=M_{u\xi,v}(u';g), 
  \end{align*}
  by \thmref{Lem}{hcmod-pairid}. This proves the assertion.
\end{proof}

We now again take up the proof of our main theorem.

\begin{proof}[\prfof{Th}{cw} (continued)]
  For $v'\in E$, we have 
  \[
    v'=0\ \Longleftrightarrow\ \forall\xi\in\tilde V\,:\,\Dual0\xi{v'}=0.
  \]
  Hence, by \thmref{Lem}{matrixcoeff-l-id}, the subspace $\mathscr N$ is invariant under $L^n$, where $L$ is the regular $G$-representation. That it is invariant under left convolution by $\Sw0G$ now follows from the identity 
  \[
    R_{S(u)}(\omega)*\varpi=\omega*(L_u(\varpi))
  \]
  valid for $u\in\Uenv0{\ger g}$, $\omega\in\Sw0{G_0}$, and $\varpi\in\Sw0G$, together with Equation \eqref{eq:sw-r-iso}.

  Since $\Sw0G\cong\Sw0{G_0}\otimes\bigwedge\ger g_\odd$ is the locally convex direct sum of finitely many copies of $\Sw0{G_0}$, it follows directly from the definition in Equation \eqref{eq:swact-def} that 
  \[
    \phi:\Sw0G^n\longrightarrow E:(\omega_1,\dotsc,\omega_n)\longmapsto\sum\nolimits_j\Pi(\omega_j)v_j
  \]
  is continuous, so that $\mathscr N$ is also closed, as claimed. 

  Hence, if we define 
  \[
    V_+\defi\Sw0G^n/\mathscr N,
  \]
  then this is a continuous non-degenerate Fr\'echet $\Sw0G$-module. By \thmref{Prop}{sw-dm}, the $\Sw0G$-action is integrated from a unique $SF$-representation $\pi$ of $G$. 

  The map induced by $\phi$ identifies $V_+$ (as a super-vector space) with the subspace
  \[
    U_+\defi\Pi(\Sw0G)V=\Pi(\Sw0{G_0})V
  \]
  of $E$. By construction \cite{bk}*{\S~6}, $U_+$ is, with the quotient topology defined by the natural map $\Sw0{G_0}^n\longrightarrow U_+$ induced by $\phi$, the minimal globalisation of the module $U\in\HC(\ger g_\ev,K_0)$. But by the Casselman--Wallach theorem \cite{bk}*{Theorem 10.6}, it holds that $U_+=E_\infty$ as locally convex spaces. 

  Since $\Uenv0{\ger g}$ is $\Ad(K_0)$-locally finite, the space of $K_0\times K_0$-finite vectors is
  \[
    \Sw0G^{(K_0\times K_0)}=\Sw0{G_0}^{(K_0\times K_0)}\otimes_{\Uenv0{\ger g_\ev}}\Uenv0{\ger g}.
  \]
  From this, it is easy to deduce that $V_+$ is an $SF$-globalisation of $V$. In particular, $(V_+)|_{G_0}$ is an $SF$-globalisation of $U$. From the Casselman--Wallach theorem \cite{bk}*{Theorem 10.6} again, it follows that the map $V_+\longrightarrow U_+$ induced by $\phi$ is an isomorphism of locally convex vector spaces. In particular, $V_+$ is the space of smooth vectors of a continuous Hilbert $G$-representation.

  Now, let $F$ be any $SF$-globalisation of $V$, so that we are given an isomorphism $\psi:V\longrightarrow F^{(K_0)}$ of $(\ger g,K_0)$-modules. Invoking the Casselman--Wallach theorem (\loccit), there is a unique isomorphism $\tilde\psi:V_+\longrightarrow F$ of $SF$-representations of $G_0$ extending $\psi$. For any $u\in\Uenv0{\ger g}$, the action by $u$ on $V_+$ and $F$ is continuous. Hence, by the density of $V$ in $V_+$, it follows that $\tilde\psi$ is $\ger g$-equivariant. This shows that $V_+$ is a minimal $SF$-globalisation. The same argument shows that it is maximal, and hence follows the claim. 
\end{proof}

\section{\texorpdfstring{Application: Gel\cprime fand--Kazhdan representations}{Application: Gel'fand--Kazhdan representations}}\label{sec:gk}

In this section, we show, by way of application of our results in Section \ref{sec:cw}, that the Gel\cprime fand--Kazhdan criterion for multiplicity freeness carries over to the case of Lie supergroups. Therein, we build on the work of Sun--Zhu \cite{sun-zhu} who have shown how to present this within the framework of Lie group Casselman--Wallach theory. Antecedents are the classical results of Gel\cprime fand--Kazhdan \cite{gk} and Shalika \cite{sha}, as well as theorems of Kostant \cite{kos-wh}, Yamashita \cite{yam}, and Prasad \cite{pra}.

We retain our assumptions on the Lie supergroup $G$ from Section \ref{sec:cw}.

\begin{Def}[contra][contragredient pairs]
  A pair $(E,F)$ of continuous $G$-re\-pre\-sen\-ta\-tions is called \Define{contragredient} if there exists a $G_0$-invariant continuous bilinear map
  \[
    \Dual0\cdot\cdot:E\times F\longrightarrow\knums
  \]
  that is a perfect pairing whose restriction to $E_\infty\times F_\infty$ is $G$-invariant.

  Here, by a \Define{perfect pairing} we mean that the canonical maps
  \[
    E\longrightarrow F',\quad F\longrightarrow E'
  \]
  are isomorphisms of topological vector spaces. 
\end{Def}

\begin{Rem}
  Assume $(E_\infty,F_\infty)$ is a pair of $SF$-representations of $G$ and 
  \[
    \Dual0\cdot\cdot:E_\infty\times F_\infty\longrightarrow\knums
  \]
  is a non-degenerate continuous bilinear form that is $G$-invariant. If $U$ is a Hilbert globalisation of $E_\infty$ (which exists if $E_\infty$ is $CW$), then the space of $G_0$-smooth vectors in $F\defi E'$ coincides with $F_\infty$. Thus, $(E,F)$ a contragredient pair with underlying $SF$-representations $E_\infty$ and $F_\infty$.
\end{Rem}

As we shall presently see, contragredient pairs of representations allow for an abstract matrix coefficient map. To state this precisely, we introduce the following definition. 

\begin{Def}[temp][tempered superfunctions]
  A superfunction $f\in\Gamma(\sh O_G)$ is called \Define{tempered} if \fa $u,v\in\Uenv0{\ger g}$
  \[
    t_{u,v,N}(f)\defi\sup_{g\in G_0}s(g)^{-N}\Abs1{(L_uR_vf)(g)}<\infty
  \]
  for some $N\sge0$. Here, $s$ denotes the maximal scale, see \cite{bk}*{2.1.1}.

  The space of tempered superfunctions is denoted by $\sh T(G)$. It is topologised as the locally convex inductive limit of the spaces $\sh T_N(G)\defi\bigcap_{u,v}\{t_{u,v,N}<\infty\}$, endowed with the locally convex topology generated by the seminorms $t_{u,v,N}$, $u,v\in\Uenv0{\ger g}$.
\end{Def}

For any $\omega\in\Sw0G$, the Berezin integral 
\[
  \vphi\longmapsto\Dual0\omega\vphi\defi\int_G\omega\vphi
\]
extends uniquely to a continuous functional on $\sh T(G)$. This is easy to deduce from \thmref{Prop}{invber} and the corresponding classical facts. 

Define $\TDi0G$, the space of \Define{tempered generalised functions}, to be the strong dual of $\Sw0G$. There is a continuous linear injection
\begin{center}
  \begin{tikzcd}
    \sh T(G)\rar{}&\TDi0G.
  \end{tikzcd}
\end{center}

The following proposition generalises \cite{sun-zhu}*{Theorem 2.1}.

\begin{Prop}[matrix]
  Let $(E,F)$ be a contragredient pair of continuous $F$-re\-pre\-sen\-ta\-tions of $G$. Then the map 
  \[
    M:E_\infty\times F_\infty\longrightarrow\sh T(G),\quad M_{v\otimes v'}(u;g)\defi(-1)^{\Abs0u(\Abs0v+\Abs0v')}\Dual1{\pi_{E,0}(g)d\pi_E(u)v}{v'}
  \]
  extends continuously to a $G\times G$-equivariant separately continuous bilinear map
  \[
    M^{-\infty}:E_{-\infty}\times F_{-\infty}\longrightarrow\TDi0G,
  \]
  where $E_{-\infty}\defi (F_\infty)'$, $F_{-\infty}\defi(E_\infty)'$, and $(d\pi_E,\pi_{E,0})$ is the $G$-action on $E$.

  If, moreover, $E_\infty$ and $F_\infty$ are $CW$ $G$-representations, then $M^{-\infty}$ is continuous and the induced $(G\times G)$-equivariant continuous linear map
  \[
    E_{-\infty}\mathop{\widehat{\otimes}}\nolimits_\pi F_{-\infty}\longrightarrow\TDi0G
  \]
  is a topological morphism with closed image. 
\end{Prop}

The structure of the \emph{proof} is manifestly the same as the one given by Sun--Zhu \cite{sun-zhu}, so we shall be brief. We begin with the following lemma. 

\begin{Lem}[smoothing]
  Let $E$ be an $F$-representation of $G$. Then the bilinear map
  \[
    \Phi_E:\Sw0G\times E\longrightarrow E_\infty:(\omega,v)\longmapsto\Pi_E(\omega)v
  \]
  is well-defined and continuous. 
\end{Lem}

\begin{proof}
  That the map is well-defined follows from \thmref{Prop}{sw-dm} \eqref{swdm-i}. The continuity is an immediate consequence of \thmref{Prop}{sw-iso} and \cite{sun-zhu}*{Lemma 3.3}.
\end{proof}

\begin{proof}[\prfof{Prop}{matrix}]
  For $\omega\in\Sw0G$, we may define $\Pi^{-\infty}_E(\omega):E^{-\infty}\longrightarrow E$ by 
  \[
    \Dual1{\Pi^{-\infty}_E(\omega)v}{v'}\defi(-1)^{\Abs0\omega\Abs0v}\Dual1v{\Phi_F(\check\omega,v')},\quad v\in E^{-\infty},v'\in F.
  \]
  Then $\Pi^{-\infty}_E(\omega)$ is continuous, and the bilinear map
  \[
    \Phi^{-\infty}_E:\Sw0G\times E^{-\infty}\longrightarrow E:(\omega,v)\longmapsto\Pi^{-\infty}_E(\omega)v
  \]
  is separately continuous, both by \thmref{Lem}{smoothing}. Applying \thmref{Prop}{sw-iso} and \cite{sun-zhu}*{Lemma 3.5}, we see that it takes values in $E_\infty$ and is separately continuous with respect to the natural topology on this space. 

  We compute for $v\in E$ and $v'\in F$: 
  \[
    \Dual1{\Phi_E^{-\infty}(\omega,v)}{v'}=(-1)^{\Abs0\omega\Abs0v}\Dual1v{\Pi_F(\check\omega)v'}=\Dual1{\Pi_E(\omega)v}{v'}=\Dual1{\Phi_E(\omega,v)}{v'},
  \]
  since for $\omega=u\otimes\varpi$, $u\in\Uenv0{\ger g}$, $\varpi\in\Sw0{G_0}$, we have 
  \begin{align*}
    \Dual1v{\Pi_F(\check\omega)v'}&=\Dual1v{\Pi_{F,0}(\check\varpi)d\pi_F(S(u))v'}\\
    &=(-1)^{\Abs0u\Abs0v}\Dual1{d\pi_E(u)\Pi_{E,0}(\varpi)v}{v'}
    =(-1)^{\Abs0\omega\Abs0v}\Dual1{\Pi_E(\omega)v}{v'},  
  \end{align*}
  in view of \thmref{Prop}{sw-dm} \eqref{swdm-i}. Thus, $\Phi_E^{-\infty}$ extends $\Phi_E$.

  Altogether, the map $M^{-\infty}:E_{-\infty}\times F_{-\infty}\longrightarrow\TDi0G$,
  \begin{equation}\label{eq:minfty-def}
    \Dual1\omega{M^{-\infty}(v,v')}\defi\Dual1{\Phi_E^{-\infty}(\omega,v)}{v'}=(-1)^{\Abs0\omega\Abs0v}\Dual1v{\Phi_F^{-\infty}(\check\omega,v')},    
  \end{equation}
  is well-defined, separately continuous, and extends $M$.

  Now, assume that $E_\infty$ and $F_\infty$ are $CW$ $G$-representations. As such, they are nuclear Fr\'echet spaces \cite{bk}*{Corollary 5.6} and hence reflexive \cite{treves}*{Corollary 3 to Proposition 50.2, Corollary to Proposition 36.9}. The same holds for $\Sw0G$, by \thmref{Prop}{sw-iso}. Thus, $E^{-\infty}$, $F^{-\infty}$, and $\TDi0G$ are strong duals of reflexive Fr\'echet spaces, and $M^{-\infty}$ is automatically continuous (\opcit, Theorem 41.1). The final statement now follows from \thmref{Cor}{topmor}.
\end{proof}

We now generalise Sun--Zhu's version of the Gel\cprime fand--Kazhdan criterion \cite{sun-zhu}*{Theorem 2.3 (i)} to Lie supergroups. 

\begin{Def}[irr][irreducible representations]
  Let $U$ be an $SF$-representation of $G$. We say that $U$ is \Define{irreducible} if there is no non-zero proper closed subspace of $U$ that is $G$-invariant. 
\end{Def}

\begin{Th}[super-gk][super Gel\cprime fand--Kazhdan criterion]
  Let $H_1,H_2$ be closed subsupergroups of $G$, $\chi_i:H_i\longrightarrow\knums^\times$ characters of $H_i$, and $\sigma:G\longrightarrow G$ an anti-automorphism. Assume that any $T\in\TDi0G_\ev$, which is at once $(H_1\times H_2)$-relatively invariant for the character $\chi_1^{-1}\otimes\chi_2^{-1}$ and a joint eigenvector of all $D\in\Uenv0{\ger g}^G_\ev$, is fixed by $\sigma$. 

  Then, for any contragredient pair $(E,F)$ of $F$-representions of $G$ \scth $E_\infty,F_\infty$ are \emph{irreducible} $CW$ $G$-representions, we have
  \[
    \dim\Hom[_{H_1}]0{E_\infty,\chi_1}\dim\Hom[_{H_2}]0{F_\infty,\chi_2}\sle1.
  \]
  Here, $\mathrm{Hom}_H$ denotes continuous even linear maps that are equivariant with respect to the supergroup $H$.
\end{Th}

\begin{proof}
  Again, our argument is largely that of Sun--Zhu \cite{sun-zhu}, with appropriate modifications and references to our results. Let 
  \[
    0\neq v\in\GHom[_{H_1}]0{E_\infty,\chi_1}\subseteq F_{-\infty},\quad
    0\neq u\in\GHom[_{H_1}]0{F_\infty,\chi_2}\subseteq E_{-\infty},
  \]
  and set $T\defi M^{-\infty}_{u\otimes v}\in\TDi0G$, appealing to \thmref{Prop}{matrix}. (Here, $\underline{\mathrm{Hom}}_{H_1}$ denotes the space of $H_1$-equivariant continuous linear maps.)

  For $D\in\Uenv0{\ger g}$ and $\omega\in\Sw0G$, we compute
  \begin{align*}
    \Dual1\omega{DT}&=(-1)^{\Abs0D\Abs0\omega}\Dual1{R_D\omega}{M^{-\infty}_{u\otimes v}}\\
    &=(-1)^{\Abs0D\Abs0\omega}\Dual1{\Pi_E^{-\infty}(R_D\omega)u}v
    =\Dual1{\Pi_E^{-\infty}(\omega)d\pi_E^{-\infty}(D)u}v,   
  \end{align*}
  by the use of Equations \eqref{eq:minfty-def} and \eqref{eq:sw-convdef}. If now $D$ is even and $G$-invariant, then $D$ commutes with the $G$-action on $E^{-\infty}$. 

  The Harish-Chandra $(\ger g,K_0)$-module $E_\infty^{(K_0)}$ is countable-dimensional, and $\Uenv0{\ger g}$ acts irreducibly, hence Dixmier's Lemma \cite{wallach1}*{Lemma 0.5.2} applies, and $S(D)$ acts by a scalar. Since $\smash{E_\infty^{(K_0)}\subseteq E_\infty}$ is dense, it follows that $D$ acts by a scalar on $E_{-\infty}$. Thus, by the computation above, $T$ is an eigenvector of $D$.

  On the other hand, as a similar computation shows, $T$ is also relatively $(\chi_1^{-1}\otimes\chi_2^{-1})$-invariant under $(H_1\times H_2)$. By assumption, $T$ is fixed by $\sigma$.

  Let $\omega\in\Sw0G$ and $g\in G_0$. We compute
  \begin{align*}
    \Dual1{\Pi_E^{-\infty}(\omega)u}{\pi_F(\sigma(g))^{-\infty}v}&=\Dual1\omega{R_{\sigma(g)}T}=\Dual1\omega{R_{\sigma(T)}\sigma(T)}\\
    &=\Dual1\omega{\sigma(L_gT)}=(-1)^{\Abs0u\Abs0\omega}\Dual1{\pi_E^{-\infty}(g)u}{\Pi_F^{-\infty}(\check\omega)v}.    
  \end{align*}
  By the irreducibility of $E_\infty$ and $F_\infty$, we conclude that 
  \[
    \Pi_E^{-\infty}(\omega)u=0\ \Longleftrightarrow\ \Pi_F^{-\infty}(\check\omega)v=0.
  \]
  Hence, for any other $0\neq u'\in\GHom[_{H_1}]0{F_\infty,\chi_2}$, the continuous linear maps 
  \[
    \Sw0G\longrightarrow E_\infty:\omega\longmapsto\Pi_E^{-\infty}(\omega)u,\quad
    \Sw0G\longrightarrow E_\infty:\omega\longmapsto\Pi_E^{-\infty}(\omega)u',
  \]
  have the same kernel $W$ (say), and induce continuous linear maps
  \[
    \vphi,\vphi':\Sw0G/W\longrightarrow E_\infty.
  \]
  These are $G$-equivariant by their definition, so they are isomorphisms with closed image, by the token of \thmref{Cor}{topmor}. They are non-zero, and therefore surjective, by the assumption of irreducibility. 

  Hence, $\psi\defi\vphi'^{-1}\circ\vphi$ is a well-defined continuous even linear and $G$-equivariant automorphism of $E_\infty$. Restricted to $\smash{E_\infty^{(K_0)}}$, it is a constant, by Dixmier's Lemma (\loccit) again. This shows that $u'\in\knums u$, by applying \thmref{Lem}{approxid}. A similar argument applies to $v$, proving the assertion.
\end{proof}

\end{document}